\documentclass[reqno,11pt]{article}
\usepackage{dsfont}

\usepackage[T1]{fontenc}      
\usepackage{ae}                      
\usepackage[latin1]{inputenc}
\usepackage{lmodern}
\usepackage{amsmath}
\usepackage{amsfonts}
\usepackage[french,english]{babel}
\usepackage{amssymb}
\usepackage{dsfont}
\usepackage{array}
\usepackage{color}
\usepackage{epsfig}
\usepackage{fancybox}
\usepackage{amsthm}
\usepackage{frcursive}
\usepackage{float}
\usepackage{shorttoc}
\usepackage{color}
\usepackage[colorlinks=true,citecolor=blue]{hyperref}

\usepackage{shorttoc}
\usepackage[boxruled]{algorithm2e}
\usepackage{pdfsync}

\usepackage[colorlinks=true,citecolor=blue]{hyperref}

\RequirePackage{amsthm, amsmath, amsfonts, amssymb}%

\DeclareMathOperator{\pen}{pen}
\DeclareMathOperator{\diag}{diag}
\DeclareMathOperator{\RE}{RE}
\DeclareMathOperator{\Card}{Card}
\newcommand{\bs}{\boldsymbol}
\newcommand{\sumin}{\sum_{i=1}^n}
\newcommand{\sumip}{\sum_{j=1}^M}
\newcommand{\intt}{\int_{0}^{t}}
\newcommand{\inttau}{\int_{0}^{\tau}}

\newcommand{\e}{{\mathrm e}}
\newcommand{\diff}{\mathrm d}

\theoremstyle{plain} 
\newtheorem{theorem}{Theorem}[section]
\newtheorem{lemma}[theorem]{Lemma} 
\newtheorem{definition}[theorem]{Estimation procedure}
\newtheorem{corollary}[theorem]{Corollary} 
\newtheorem{proposition}[theorem]{Proposition} 
{\bf}{\rm}%
\newtheorem{remark}[theorem]{Remark}

\usepackage[toc,page]{appendix} 

\newcounter{hypc}
\newcommand{\hyp}[1]{\stepcounter{hypc}\tag{$\mathbf{A_{\thehypc}}$}
\label{#1}}

\newcounter{RE}
\newcommand{\hypRE}[1]{\stepcounter{RE}\tag{$\mathbf{RE(s,a_0)_{}}$}
\label{#1}}

\newcounter{REinconnu}
\newcommand{\hypREinconnu}[1]{\stepcounter{REinconnu}\tag{$\mathbf{\widetilde{RE}(s,r_0)_{}}$}
\label{#1}}

\begin{document}
\title{Oracle inequalities for the Lasso in the high-dimensional Aalen multiplicative intensity model}
\author{Sarah Lemler\\
\small{
Laboratoire Statistique et G\'enome UMR CNRS 8071- USC INRA,}\\
\small{Universit\'e d'\'Evry Val d'Essonne, France}\\
\small{\textit{e-mail :} \texttt{sarah.lemler@genopole.cnrs.fr}}
}
\date{}
\maketitle

\begin{abstract} 
In a general counting process setting, we consider the problem of obtaining a prognostic on the survival time adjusted on covariates in high-dimension. Towards this end, we construct an estimator of the whole conditional intensity. We estimate it by the best Cox proportional hazards model given two dictionaries of functions. The first dictionary is used to construct an approximation of the logarithm of the baseline hazard function and the second to approximate the relative risk. We introduce a new data-driven weighted Lasso procedure to estimate the unknown parameters of the best Cox model approximating the intensity. We provide non-asymptotic oracle inequalities for our procedure in terms of an appropriate empirical Kullback divergence. Our results rely on an empirical Bernstein's inequality for martingales with jumps and properties of modified self-concordant functions. 

\begin{flushleft}
\end{flushleft}

\end{abstract}
%

\begin{small}
\tableofcontents
\end{small}
%
%
\section{Introduction}
We consider one of the statistical challenges brought by the recent advances in biomedical technology to clinical applications. For example, in Dave et al.~\cite{Dave}, the considered data relate 191 patients with follicular lymphoma. The observed variables are the survival time, that can be right-censored, clinical variables, as the age or the disease stage, and 44 929 levels of gene expression. In this high-dimensional right-censored setting, there are two clinical questions. One is to determine prognostic biomarkers, the second is to predict the survival from follicular lymphoma adjusted on covariates. We focus our interest on the second (see Gourlay \cite{Gourlay} and Steyerberg \cite{Steyerberg}). 
As a consequence, we consider the statistical question of estimating the whole conditional intensity.
To adjust on covariates, the most popular semi-parametric regression model is the Cox proportional hazards model (see Cox \cite{Cox}) : the conditional hazard rate function of the survival time $T$ given the vector of covariates $\bs{Z}=(Z_{1},...,Z_{p})^T$ is defined by
\begin{align}
\label{eq:Cox}
\lambda_{0}(t,\bs{Z})=\alpha_{0}(t)\exp(\bs{\beta_{0}^TZ}),
\end{align}
where $\bs{\beta_{0}}=(\beta_{{0}_{1}},...,\beta_{{0}_{p}})^T$ is the vector of regression coefficients and $\alpha_{0}$ is the baseline hazard function. The unknown parameters of the model are $\bs{\beta_{0}}\in\mathbb{R}^p$ and the function $\alpha_{0}$. To construct an estimator of $\lambda_0$, one usually considers the partial likelihood introduced by Cox \cite{Cox} to derive an estimator of $\bs{\beta_0}$ and then plug this estimator to obtain the well-known Breslow estimator of $\alpha_0$. We propose in this paper an alternative one-step strategy. 


\subsection{Framework}
Before describing our strategy, let us clarify our framework.
We consider the general setting of counting processes. For $i=1,...,n$, let $N_{i}$ be a marked counting process and $Y_{i}$ a predictable random process with values in $[0,1]$. Let $(\Omega,\mathcal{F},\mathbb{P})$ be a probability space and $(\mathcal{F}_{t})_{ t\geq 0}$ be the filtration 
defined by 
$$\mathcal{F}_{t}=\sigma\{N_{i}(s),Y_{i}(s), 0\leq s\leq t, \bs{Z_{i}}, i=1,...,n\},$$
where $\bs{Z_i}=(Z_{i,1},...,Z_{i,p})^T\in\mathbb{R}^p$ is the $\mathcal{F}_{0}$-measurable random vector of covariates of individual $i$. 
Let $\Lambda_{i}(t)$ be the compensator of the process $N_{i}(t)$ with respect to $(\mathcal{F}_{t})_{t\geq 0}$, so that $M_{i}(t)=N_{i}(t)-\Lambda_{i}(t)$ is a $(\mathcal{F}_{t})_{t\geq 0}$-martingale. 

%
%

\begin{align}
&  \text{The process }N_{i} \text{ satisfies the Aalen multiplicative intensity model : for all } t\geq 0, \notag\\
&\hyp{eq:lambda} \hspace{4cm}\Lambda_{i}(t)=\displaystyle{\int_{0}^{t}\lambda_{0}(s,\bs{Z_{i}})Y_{i}(s)\diff s}, \notag \\
& \mbox{where } \lambda_{0} \mbox{ is an unknown nonnegative function called intensity. } \notag
\end{align}


This general setting, introduced by Aalen \cite{aalen}, embeds several particular examples as censored data, marked Poisson processes and Markov processes (see Andersen et al. \cite{ABG} for further details). 
\begin{remark}
\label{censoring}
In the specific case of right censoring, let $(T_{i})_{i=1,...,n}$ be i.i.d. survival times of $n$ individuals and $(C_{i})_{i=1,...,n}$ their i.i.d. censoring times. We observe $\{(X_{i},\bs{Z_{i}},\delta_{i})\}_{i=1,...,n}$ where $X_{i}=\min({T_{i}},{C_{i}})$ is the event time, $\bs{Z_{i}}=(Z_{i,1},...,Z_{i,p})^T$ is the vector of covariates and $\delta_{i}=\mathds{1}_{\{T_{i}\leq C_{i}\}}$ is the censoring indicator. The survival times $T_{i}$ are supposed to be conditionally independent of the censoring times $C_{i}$ given some vector of covariates $\bs{Z_{i}}=(Z_{i,1},...,Z_{i,p})^T\in\mathbb{R}^p$ for ${i=1,...,n}$. With these notations, the $(\mathcal{F}_{t})$-adapted processes $Y_{i}$ and $N_{i}$ are respectively defined as the at-risk process $Y_{i}(t)=\mathds{1}_{\{X_{i}\geq t\}}$ and the counting process $N_{i} (t)=\mathds{1}_{\{X_{i}\leq t,\delta_{i}=1\}}$ which jumps when the ith individual dies. 
\end{remark}

We observe the independent and identically distributed (i.i.d.) data $(\bs{Z_{i}}, N_{i}(t), Y_{i}(t),
 i=1,...,n, 0\leq t\leq \tau)$, where $[0,\tau]$ is the time interval between the beginning and the end of the study. 
 
 \begin{align}
& \hyp{A_0} \text{On }[0,\tau], \text{ we assume that }  A_0=\underset{1\leq i\leq n}{\sup}\Big\{\displaystyle{\int_{0}^{\tau}\lambda_{0}(s,\bs{Z_{i}})\diff s}\Big\}<\infty. \notag 
\end{align}
This is the standard assumption in statistical estimation of intensities of counting processes, see Andersen et al. \cite{ABG} for instance. We also precise that, in the following, we work conditionally to the covariates and from now on, all probabilities $\mathbb{P}$ and expectations $\mathbb{E}$ are conditional to the covariates.
Our goal is to estimate $\lambda_0$ non-parametrically in a high-dimensional setting, i.e. when the number of covariates $p$ is larger than the sample size $n$ ($p\gg n$). 

\subsection{Previous results}

In high-dimensional regression, the benchmarks for results are the ones obtained in the additive regression model. In this setting, Tibshirani \cite{Tib} has introduced the Lasso procedure, which consists in minimizing an $\ell_{1}$-penalized criterion.
%
The Lasso estimator has been widely studied for this model, with consistency results (see Meinshausen and B\"uhlmann \cite{Meinshausen}) and variable selection results (see Zhao and Yu \cite{Zhao}, Zhang and Huang \cite{Zhang1}). Recently, attention has been directed on establishing non-asymptotic oracle inequalities for the Lasso (see Bunea et al. \cite{BTW1,BTW2}, Bickel et al. \cite{BRT}, Massart and Meynet \cite{MM}, Bartlett \cite{Mendelson12} and Koltchinskii \cite{Koltchinskii11} among others). 

In the setting of survival analysis, the Lasso procedure has been first considered by Tibshirani \cite{T}  and applied to the partial log-likelihood. 
More generally, other procedures have been introduced for the parametric part of the Cox model : the adaptive Lasso, the smooth clipped absolute deviation penalizations and the Danzig selector are respectively considered in Zou \cite{Zou08}, Zhang and Lu \cite{ZhangLu07}, Fan and Li \cite{FanLi02} and Antoniadis et al. \cite{AFL}. Non parametric approaches are considered in Letu\'e \cite{Letue00}, Hansen et al. \cite{Hansen} and Comte et al. \cite{CGG}. Lasso procedures for the alternative Aalen additive model have been introduced in Martinussen and Scheike \cite{Martinussen09} and Ga\"iffas and Guilloux \cite{Agathe}.

All of the existing results in the Cox model are based on the partial log-likelihood, which does not answer the clinical question associated to a prognosis. Antoniadis et al. \cite{AFL} have established asymptotic estimation inequalities in the Cox proportional hazard model for the Dantzig estimator (see Bickel et al. \cite{BRT} for a comparison between these two estimators in an additive regression model). In Bradic et al. \cite{Fan}, asymptotic estimation inequalities for the Lasso estimator have also been obtained in the Cox model. 
More recently, Kong and Nan \cite{Kong} and Bradic and Song \cite{bradic2012} have established non-asymptotic oracle inequalities for the Lasso in the generalized Cox model
\begin{equation}
\label{eq:CoxGen}
\lambda_0(t,\bs{Z})=\alpha_0(t)\exp(f_0(\bs{Z})),
\end{equation}
where $\alpha_0$ is the baseline hazard function and $f_0$ a function of the covariates. However, the focus in both papers is on the Cox partial log-likelihood, the obtained results are either on $f_{\hat\beta_{L}}-f_{0}$ or on $\bs{\hat\beta_{L}}-\bs{\beta_{0}}$ for $f_0(\bs{Z})=\bs{\beta_0^TZ}$  and the problem of estimating the whole intensity $\lambda_{0}$ is not considered, as needed for the prevision of the survival time. 

\subsection{Our contribution}
The first motivation of the present paper is to address the problem of estimating $\lambda_{0}$  defined in (\ref{eq:lambda}) regardless of an underlying model. We use an agnostic learning approach, see Kearns et al. \cite{Kearns94}, to construct an estimator that mimics the performance of the best Cox model, whether this model is true or not. More precisely, we will consider candidates for the estimation of $\lambda_0$ of the form 
\begin{equation*}
\lambda_{\beta,\gamma}(t,\bs{Z})=\alpha_{\gamma}(t)\e^{f_{\beta}(\bs{Z})} \text{ for } (\bs\beta,\bs\gamma)\in\mathbb{R}^M\times\mathbb{R}^N,
\end{equation*}
where $f_\beta$ and $\alpha_\gamma$ are respectively linear combinations of functions of two dictionaries $\mathbb{F}_{M}$ and $\mathbb{G}_{N}$. The estimator of $\lambda_0$ is defined as the candidate which minimizes a weighted $\ell_1$-penalized total log-likelihood as opposed to the Cox partial log-likelihood. 
The second motivation of the paper is to obtain non-asymptotic oracle inequalities for Lasso estimators of the complete intensity $\lambda_0$. Indeed, in practice, one can not consider that the asymptotic regime has been reached, cf. in Dave et al. \cite{Dave} for example. In addition, Comte et al. \cite{CGG} established non-asymptotic oracle inequalities for the whole intensity but not in a high-dimensional setting and to the best of our knowledge, no non-asymptotic results for the estimation of the whole intensity in high dimension exist in the literature.

Towards this end, we will proceed in two steps. In a first step, we assume that $\lambda_{0}$ verifies Model (\ref{eq:CoxGen}), where $\alpha_{0}$ is assumed to be known. In this particular case, the only nonparametric function to estimate is $f_{0}$ and we estimate it by a linear combination of functions of the dictionary $\mathbb{F}_{M}$. In this setting, we obtain non-asymptotic oracle inequalities for the Cox model when $\alpha_{0}$ is supposed to be known. In a second step, we consider the general problem of estimating the whole intensity $\lambda_{0}$. We state non-asymptotic oracle inequalities both in terms of empirical Kullback divergence and weighted empirical quadratic norm for our Lasso estimators, thanks to properties of modified self-concordant functions (see Bach \cite{Bach}).

These results are obtained via three ingredients : a new Bernstein's inequality, a modified Restricted Eigenvalue condition on the expectation of the weighted Gram matrix and modified self-concordant functions. Let us be more precise. We establish empirical versions of Bernstein's inequality involving the optional variation for martingales with jumps (see Ga\"iffas and Guilloux \cite{Agathe} and Hansen et al. \cite{Hansen} for related results). This allows us to define a fully data-driven weighted $\ell_{1}$-penalization. 
For the resulting estimator, we work under a modified Restricted Eigenvalue condition according to which the expectation of a weighted Gram matrix fullfilled the Restricted Eigenvalue condition (see Bickel et al. \cite{BRT}). This new version of the Restricted Eigenvalue condition is both new and weaker than the comparable condition in the Cox model. Finally, we extend the notion of self-concordance (see Bach \cite{Bach}) to the problem at hands in order to connect our weighted empirical quadratic norm and our empirical Kullback divergence. In this context, we state the first fast non-asymptotic oracle inequality for the whole intensity.

The paper is organized as follows. In Section 2, we describe the framework and the Lasso procedure for estimating the intensity. The estimation risk that we consider and its associated loss function are presented. In Section 3, prediction and estimation oracle inequalities in the particular Cox model with known baseline hazard function are stated. In Section 4, non-asymptotic oracle inequalities with different convergence rates are given for a general intensity. Section 5 is devoted to statement of empirical Bernstein's inequalities associated to our processes. Proofs are gathered in section 6.

\section{Estimation procedure}
\subsection{The estimation criterion and the loss function}

To estimate the intensity $\lambda_0$, we consider the total empirical log-likelihood. 
By Jacod's Formula (see Andersen et al. \cite{ABG}), the log-likelihood based on the data $(\bs{Z_{i}}, N_{i}(t), Y_{i}(t), i=1,...,n, 0\leq t\leq \tau)$ is given by  
\begin{equation*}
C_{n}(\lambda)=-\dfrac{1}{n}\displaystyle{\sum_{i=1}^{n}\left\{\int_{0}^{\tau}\log\lambda(t,\bs{Z_{i}})\diff N_{i}(t)-\int_{0}^{\tau}\lambda(t,\bs{Z_{i}})Y_{i}(t)\diff t\right\}}.
\end{equation*}
Our estimation procedure is based on the minimization of this empirical risk.
To this empirical risk, we associate the empirical Kullback divergence defined by 
\begin{align}
\nonumber
\widetilde{K}_{n}(\lambda_{0},\lambda)&=\dfrac{1}{n}\displaystyle{\sum_{i=1}^{n}\int_{0}^{\tau}\left(\log\lambda_{0}(t,\bs{Z_{i}})-\log\lambda(t,\bs{Z_{i}})\right)\lambda_{0}(t,\bs{Z_{i}})Y_{i}(t)\diff t}\\
\label{eq:Kullback}
&-\dfrac{1}{n}\displaystyle{\sum_{i=1}^{n}\int_{0}^{\tau}\left(\lambda_{0}(t,\bs{Z_{i}})-\lambda(t,\bs{Z_{i}})\right)Y_{i}(t)\diff t}.
\end{align}
We refer to van de Geer \cite{SVG} and Senoussi \cite{Senoussi} for close definitions. Notice in addition, that this loss function is closed to the Kullback-Leibler information considered in the density framework (see Stone \cite{Stone94} and Cohen and Le Pennec \cite{LePennecCohen13}).
The following proposition justify the choice of this criterion.
\begin{proposition}
\label{prop:Kullback}
The empirical Kullback divergence $\widetilde{K}_{n}(\lambda_{0},\lambda)$ is nonnegative and equals zero if and only if $\lambda=\lambda_{0}$ almost surely on the interval $[0,\tau\wedge\sup\{t : \exists i \in\{1,...,n\}, Y_i(t)\neq 0\}]$.
\end{proposition}

\begin{remark}
In the specific case of right censoring, the proposition holds true on $[0,\tau\wedge\underset{1\leq i\leq n}{\max}X_i]$. In this case, we can specify that $\mathbb{P}([0,\tau]\subset[0,\underset{1\leq i\leq n}{\max}X_i])=1-(1-S_{T}(\tau))^n(1-S_{C}(\tau))^n$, where $S_{T}$ and $S_C$ are the survival functions of the survival time $T$ and the censoring time $C$ respectively. From \ref{A_0}, $S_T(\tau)>0$ and if $\tau$ is such that $S_C(\tau)>0$, then $\mathbb{P}([0,\tau]\subset[0,\underset{1\leq i\leq n}{\max}X_i])$ is large. See Gill \cite{Gill83} for a discussion on the role of $\tau$. 
\end{remark}

In the following, we consider that we estimate $\lambda_0(t)$ for $t$ in $[0,\tau\wedge\sup\{t : \exists i \in\{1,...,n\}, Y_i(t)\neq 0\}]$. Let introduce the weighted empirical quadratic norm defined for all function $h$ on $[0,\tau]\times\mathbb{R}^p$ by
\begin{equation}
\label{eq:norm}
||h||_{n,\Lambda}=\sqrt{\dfrac{1}{n}\sumin\inttau (h(t,\bs{Z_{i}}))^2\diff \Lambda_{i}(t)},
\end{equation}
where $\Lambda_{i}$ is defined in (\ref{eq:lambda}). Notice that, in this definition, the higher the intensity of the process $N_{i}$ is, the higher the contribution of individual $i$ to the empirical norm is. This norm is connected to the empirical Kullback divergence, as it will be shown in Proposition \ref{prop:comparaison1}. Finally, for a vector $\bs{b}$ in $\mathbb{R}^M$, we define, $||\bs b||_1=\sum_{j=1}^{M}|b_j|$ and $||\bs b||_2^2=\sum_{j=1}^{M}b_j^2$.

\subsection{Weighted Lasso estimation procedure}

The estimation procedure is based on the choice of two finite sets of functions, called dictionaries.
Let $\mathbb{F}_M=\{f_{1},...,f_{M}\}$ where $f_{j} : \mathbb{R}^p \rightarrow \mathbb{R}$ for $j=1,...,M$, and $\mathbb{G}_{N}=\{\theta_{1},...,\theta_{N}\}$ where $\theta_{k} : \mathbb{R}_{+} \rightarrow \mathbb{R}$ for  $k=1,...,N$, be two dictionaries. Typically the size of the dictionary $\mathbb{F}_M$ used to estimate the function of the covariates in a high-dimensional setting is large, i.e. $M\gg n$, whereas to estimate a function on $\mathbb{R}_+$, we consider a dictionary $\mathbb{G}_N$ with size $N$ of the order of $n$. The sets $\mathbb{F}_M$ and $\mathbb{G}_{N}$ can be collections of functions such as wavelets, splines, step functions, coordinate functions etc. They can also be collections of several estimators computed using different tuning parameters. To make sure that no identification problems appear by using two dictionaries, it is assumed that only the dictionary $\mathbb{G}_{N}=\{\theta_{1},...,\theta_{N}\}$ can contain the constant function, not $\mathbb{F}_M=\{f_{1},...,f_{M}\}$.
The candidates for the estimator of $\lambda_0$ are of the form
\begin{equation*}
\lambda_{\beta,\gamma}(t,\bs{Z_{i}})=\alpha_{\gamma}(t)\e^{f_{\beta}(\bs{Z_{i}})} \text{ with } \log\alpha_{\gamma}=\displaystyle{\sum_{k=1}^{N}\gamma_{k}\theta_{k}} \text{ and } f_{\beta}=\displaystyle{\sum_{j=1}^{M}\beta_{j}f_{j}}.
\end{equation*}


The dictionaries $\mathbb{F}_{M}$ and $\mathbb{G}_{N}$ are chosen such that the two following assumptions are fullfiled. 

\begin{align}
&  \hyp{ass:dico1}\text{For all } j \text{ in } \{1,...,M\}, ||f_{j}||_{n,\infty}=\underset{1\leq i\leq n}{\max}|f_{j}(Z_{i})|<\infty. \notag 
\end{align}

\begin{align}
&  \hyp{ass:dico2}\text{For all } k \text{ in } \{1,...,N\}, ||\theta_{k}||_{\infty}=\underset{t\in[0,\tau]}{\max}|\theta_{k}(t)|<\infty. \notag 
\end{align}




We consider a weighted Lasso procedure for estimating $\lambda_0$. 
\begin{definition}
\label{EstimationProc}
The Lasso estimator of $\lambda_0$ is defined by $\lambda_{\hat\beta_{L},\hat\gamma_{L}}$, where
\begin{equation*}
(\bs{\hat{\beta}_{L}},\bs{\hat\gamma_{L}})=\underset{(\bs{\beta},\bs\gamma)\in\mathbb{R}^M\times\mathbb{R}^N}{\arg\min}\{C_{n}(\lambda_{\beta,\gamma})+\pen(\bs{\beta})+\pen(\bs\gamma)\},
\end{equation*}
with 
\begin{equation*}
\pen(\bs{\beta})=\displaystyle{\sum_{j=1}^{M}\omega_{j}|\beta_{j}}| \text{ and } \pen(\bs\gamma)=\displaystyle{\sum_{k=1}^{N}\delta_{k}|\gamma_{k}}|.
\end{equation*}
\end{definition}
The positive data-driven weights $\omega_{j}=\omega(f_j,n,M,\nu,x)$, $j=1,...,M$ and $\delta_{k}=\delta(\theta_k,n,N,\tilde\nu,y)$, $k=1,...,N$ are defined as follows. Let $x>0$, $y>0$, $\varepsilon>0$, $\tilde\varepsilon>0$, $c=2\sqrt{2(1+\varepsilon)}$, $\tilde c=2\sqrt{2(1+\tilde\varepsilon)}$ and $(\nu,\tilde\nu)\in(0,3)^2$ such that $\nu>\Phi(\nu)$ and $\tilde\nu>\Phi(\tilde\nu)$, where $\Phi(u)=\exp(u)-u-1$. With these notations, the weigths are defined by 


\begin{align}
\label{eq:weights}
\omega_{j}=c\sqrt{\dfrac{\hat{W}_n^{\nu}(f_j)(x+\log M)}{n}}+2\dfrac{x+\log M}{3n}||f_j||_{n,\infty} \quad \hspace{-0.1cm}\text{and} \hspace{-0.1cm}\quad \delta_{k}=\tilde c\sqrt{\dfrac{\hat{T}_n^{\tilde\nu}(\theta_k)(y+\log N)}{n}}+2\dfrac{y+\log N}{3n}||\theta_k||_{\infty},
\end{align}
for
\begin{align}
 \hat W_n^{\nu}(f_j)&=\dfrac{\nu/n}{\nu/n-\Phi(\nu/n)}\hat V_{n}(f_j)+\dfrac{x/n}{\nu/n-\Phi(\nu/n)}||f_j||_{n,\infty}^2,\\  
\hat T_n^{\tilde\nu}(\theta_k)&=\dfrac{\tilde\nu/n}{\tilde\nu/n-\Phi(\tilde\nu/n)}\hat R_{n}(\theta_k)+\dfrac{y/n}{\tilde\nu/n-\Phi(\tilde\nu/n)}||\theta_k||_{\infty}^2,
\end{align}
where $\hat V_{n}(f_{j})$ and $\hat R_{n}(\theta_{k})$ are the "observable" empirical variance of $f_{j}$ and $\theta_{k}$ respectively, given by
\begin{equation*}
\label{eq:obsempvar}
\hat V_{n}(f_{j})=\dfrac{1}{n}\sumin\inttau(f_{j}(\bs{Z_{i}}))^2\diff N_{i}(s) \text{ and }\hat R_{n}(\theta_{k})=\dfrac{1}{n}\sumin\inttau  (\theta_{k}(s))^2\diff N_{i}(s).
\end{equation*}
\begin{remark}
The general Lasso estimator for $\bs\beta$ is classically defined by
$$\bs{\hat\beta_{L}}=\underset{\bs\beta\in\mathbb{R}^M}{\arg\min}\{C_{n}(\lambda_{\beta})+\Gamma\sumip|\beta_{j}|\},$$
with $\Gamma>0$ a smoothing parameter. Usually, $\Gamma$ is of order $\sqrt{\log M/n}$ (see Massart and Meynet \cite{MM} for the usual additive regression model and Antoniadis et al. \cite{AFL} for the Cox model among other). The Lasso penalization for $\bs\beta$ corresponds to the simple choice $\omega_{j}=\Gamma$ where $\Gamma>0$ is a smoothing parameter. Our weights could be compared with those of Bickel and al. \cite{BRT} in the case of an additive regression model with a gaussian noise. They have considered a weighted Lasso with a penalty term of the form $\Gamma\sum_{j=1}^{M}||f_j||_n|\beta_j|$, with $\Gamma$ of order $\sqrt{\log M/n}$ and $||.||_n$ the usual empirical norm. 
We can deduce from the weights $\omega_j$ defined by (\ref{eq:weights}) higher suitable weights that can be written $\Gamma_{n,M}^1\tilde\omega_j$ with $\tilde\omega_j=\sqrt{\hat W_n^{\nu}(f_j)}$, which is of order $\sqrt{\hat V_n(f_j)}$ and
$$\Gamma_{n,M}^1= c\sqrt{\dfrac{x+\log M}{n}}+2\dfrac{x+\log M}{3n}\underset{1\leq j\leq M}{\max}\dfrac{||f_{j}||_{n,\infty}}{\sqrt{\hat{W}_n^{\nu}(f_j)}}.$$
The regularization parameter $\Gamma_{n,M}^1$ is still of order $\sqrt{\log M/n}$. The weights $\tilde\omega_j$ correspond to the estimation of the weighted empirical norm $||.||_{n,\Lambda}$ that is not observable and play the same role than the empirical norm $||f_j||_n$ in Bickel et al. \cite{BRT}. These weights are also of the same form as those of van de Geer \cite{SVG1} for the logistic model. 

The idea of adding some weights in the penalization comes from the adaptive Lasso, although it is not the same procedure. Indeed, in the adaptive Lasso (see Zou \cite{Zou}) one chooses $\omega_{j}=|\tilde\beta_{j}|^{-a}$ where $\tilde\beta_{j}$ is a preliminary estimator and $a>0$ a constant. The idea behind this is to correct the bias of the Lasso in terms of variables selection accuracy (see Zou \cite{Zou} and Zhang \cite{Zhang} for regression analysis and Zhang and Lu \cite{Lu} for the Cox model). The weights $\omega_{j}$ can also be used to scale each variable at the same level, which is suitable when some variables have a large variance compared to the others. 

\end{remark}


\section{Oracle inequalities for the Cox model when the baseline hazard function is known}
As a first step, we suppose that the intensity satisfies the generalization of the Cox model (\ref{eq:CoxGen}) with a known baseline function $\alpha_0$.
In this context, only $f_{0}$ has to be estimated and $\lambda_{0}$ is estimated~ by  
\begin{equation}
\label{eq:lassosimple}
\lambda_{\hat\beta_{L}}(t,\bs{Z_{i}})=\alpha_{0}(t)\e^{f_{\hat\beta_{L}}(\bs{Z_{i}})} \text{ and }\bs{\hat\beta_{L}}=\underset{\bs\beta\in\mathbb{R}^M}{\arg\min}\{C_{n}(\lambda_{\beta})+\pen(\bs\beta)\}.
\end{equation}
In this section, we state non-asymptotic oracle inequalities for the prediction loss of the Lasso in terms of the Kullback divergence. These inequalities allow us to compare the prediction error of the estimator and the best approximation of the regression function by a linear combination of the functions of the dictionary in a non-asymptotic way. 

\subsection{A slow oracle inequality}

In the following theorem, we state an oracle inequality in the Cox model with slow rate of convergence, i.e. with a rate of convergence of order $\sqrt{\log M/n}$. This inequality is obtained under a very light assumption on the dictionary $\mathbb{F}_{M}$.

\begin{proposition}
\label{thm:slow-oracle}
Consider Model (\ref{eq:CoxGen}) with known $\alpha_{0}$. Let $x>0$ be fixed, $\omega_{j}$ be defined by (\ref{eq:weights}) and  for $\bs{\beta}\in\mathbb{R}^M,$ 
$$\pen(\bs\beta)=\displaystyle{\sumip\omega_{j}|\beta_{j}|}.$$
Let $A_{\varepsilon,\nu}$ be some numerical positive constant depending only on $\varepsilon$ and $c_\ell$, and $x>0$ be fixed. 
Under Assumption \ref{ass:dico1}, with a probability larger than $1-A_{\varepsilon,\nu}\e^{-x}$, then 
\begin{align}
\label{eq:slow-oracle}
\widetilde{K}_{n}(\lambda_{0},\lambda_{\hat{\beta}_{L}})\leq\underset{\beta\in\mathbb{R}^M}{\inf}\left(\widetilde{K}_{n}(\lambda_{0},\lambda_{\beta})+2\pen(\bs{\beta})\right).
\end{align}
\end{proposition}

This theorem states a non-asymptotic oracle inequality in prediction on the conditional hazard rate function in the Cox model.
The $\omega_{j}$ are the order of $\sqrt{\log M/n}$ and the penalty term is of order 
$||\bs{\beta}||_{1}\sqrt{\log M/n}$. This variance order is usually referred as a slow rate of convergence in high dimension (see Bickel et al. \cite{BRT} for the additive regression model, Bertin et al. \cite{BlPR} and Bunea et al. \cite{BTWB} for density estimation).


\subsection{A fast oracle inequality}
\label{subsec:fast}

Now, we are interested in obtaining a non-asymptotic oracle inequality with a fast rate of convergence of order $\log M/n$ and we need further assumptions in order to prove such result.  
In this subsection, we shall work locally, for $\mu>0$, on the set $\Gamma_M(\mu)=\{\bs\beta\in\mathbb{R}^M : ||\log\lambda_{\bs\beta}-\log\lambda_0||_{n,\infty}\leq \mu\}$, simply denoted $\Gamma(\mu)$ to simplify the notations and we consider the following assumption :

\begin{align}
&\hyp{ass:voisinage}  \text{There exists }\mu>0, \text{ such that } \Gamma(\mu) \text{ contains a non-empty open set of } \mathbb{R}^M.  \notag
\end{align}

This assumption has already been considered by van de Geer \cite{SVG1} or Kong and Nan \cite{Kong}. Roughly speaking, it means that one can find a set where we can restrict our attention for finding good estimator of $f_0$. This assumption is needed in order to connect, via the notion of self-concordance (see Bach \cite{Bach}),  the weighted empirical quadratic norm and the empirical Kullback divergence (see Proposition \ref{prop:comparaison}).
%

The weighted Lasso estimator becomes
\begin{equation}
\label{eq:lassoGamma}
\bs{\hat\beta_L^{\mu}}=\underset{\bs\beta\in\Gamma(\mu)}{\arg\min}\{C_n(\lambda_{\bs\beta})+\pen(\bs\beta)\}.
\end{equation}
%
By definition, this weighted Lasso estimator is obtained on a ball centered around the true function $\lambda_0$. However in Assumption \ref{ass:voisinage}, we can always consider a large radius $\mu$, which weakens it. This could not change the rate of convergence in the oracle inequalities ($\sim\log M/n$) but only the range of a constant. In the particular case in which $\log\lambda_{\bs{\beta}}$ for all $\bs{\beta}\in\mathbb{R}^M$ and $\log\lambda_0$ are bounded, there exists $\mu>0$ such that $||\log \lambda_{\bs\beta}-\log\lambda_0||_{n,\infty}\leq ||\log \lambda_\beta||_{n,\infty}+||\log\lambda_0||_{n,\infty}\leq \mu$.

To achieve a fast rate of convergence, one needs an additional assumption on the Gram matrix. We choose to work under a Restricted Eigenvalue condition, as introduced in Bickel et al. \cite{BRT} for the additive regression model. This condition is one of the weakest assumption on the design matrix. See B\"uhlmann and van de Geer \cite{BvdG} and Bickel et al. \cite{BRT} for further details on assumptions required for oracle inequalities.

Let us first introduce further notations ~:
$$\bs{\Delta}=\bs{D}(\bs{\hat\beta_L^{\mu}}-\bs{\beta}) \text{ with }  \bs\beta\in\Gamma(\mu) \text{ and }\bs{D}=(\diag(\omega_{j}))_{1\leq j\leq M},$$

$$\bs{X}=(f_{j}(\bs{Z_{i}}))_{i,j}, \text{ with } i\in\{1,...,n\} \text{ and } j\in\{1,...,M\},$$

\begin{equation}
\label{Gn}
\bs{G_{n}}=\dfrac{1}{n}\bs{X^TCX} \text{ with } \bs{C}=(\diag(\Lambda_{i}(\tau)))_{1\leq i\leq n}.
\end{equation}
In the matrix $\bs{G_{n}}$, the covariates of individual $i$ is re-weighted by its cumulative risk $\Lambda_{i}(\tau)$, which is consistent with the definition of the empirical norm in (\ref{eq:norm}).
Let also $J(\bs{\beta})$ be the sparsity set  of vector  $\bs{\beta}\in\Gamma(\mu)$  defined by $J(\bs{\beta})=\{j\in\{1,...,M\} : \beta_{j}\neq0\}$, and the sparsity index is then given by $|J(\bs{\beta})|=\Card\{J(\bs\beta)\}.$
For $J\subset\{1,...,M\}$, we denote by $\bs\beta_{J}$ the vector $\bs\beta$ restricted to the set $J$ : $(\beta_{J})_{j}=\beta_{j}$ if $j\in J$ and $(\beta_{J})_{j}=0$ if $j\in J^c$ where $J^c=\{1,...,M\}\setminus J$.

Usually, in order to obtain a fast oracle inequality, we need to assume a Restricted Eigenvalue condition on the Gram matrix $\bs{G_n}$. However, since $\bs{G_n}$ is random in our case, we impose the  Restricted Eigenvalue condition to $\mathbb{E}(\bs{G_n})$, where the expectation is taken conditionally to the covariates.
\begin{align}
&  \text{For some integer } s\in\{1,...,M\} \text{ and a constant } a_{0}>0, \text{ the following condition holds :}\notag\\
&\hspace{2cm}\hypRE{ass:RE} 0<\bs\kappa_0(s,a_{0})=\underset{|J|\leq s}{\underset{J\subset\{1,...,M\},}{\min}}\underset{||\bs{b}_{J^c}||_{1}\leq a_{0}||\bs{b}_{J}||_{1}}{\underset{\bs{b}\in\mathbb{R}^M\backslash\{0\},}{\min}}\dfrac{(b^T\mathbb{E}(\bs{G_n})b)^{1/2}}{||\bs{b}_{J}||_{2}}. \notag 
\end{align}
The integer $s$ here plays the role of an upper bound on the sparsity $|J(\bs\beta)|$ of a vector of coefficients $\bs\beta$.

This assumption is weaker than the classical one and the following lemma implies that if the Restricted Eigenvalue condition is verified for $\mathbb{E}(\bs{G_n})$, then the empirical version of the Restricted Eigenvalue condition applied to $\bs{G_n}$ holds true with large probability. This modified Restricted Eigenvalue condition is new and this is the first time to our best knowledge that a fast-non asymptotic oracle inequality has been established under such a condition. 
\begin{lemma}
\label{REn}
Let  $L>0$ such that $\underset{1\leq j\leq M}{\max}\underset{1\leq i\leq n}{\max}|f_j(\bs{Z_i})|\leq L$. Under Assumptions \ref{A_0} and \ref{ass:RE}, we have \begin{equation}
\label{eq:REn}
0<\bs\kappa=\underset{|J|\leq s}{\underset{J\subset\{1,...,M\},}{\min}}\underset{||\bs{b}_{J^c}||_{1}\leq a_{0}||\bs{b}_{J}||_{1}}{\underset{\bs{b}\in\mathbb{R}^M\backslash\{0\},}{\min}}\dfrac{(\bs{b^TG_nb})^{1/2}}{||\bs{b}_{J}||_{2}} \mbox{ and } \bs\kappa=(1/\sqrt{2A_0})\bs\kappa_0(s,a_0),
\end{equation}
with probability larger than $1-\pi_n$, where 
$$\pi_n=2M^2\exp\Big[-\dfrac{n\bs\kappa^4}{2L^2(1+a_0)^2s(L^2(1+a_0)^2s+\bs\kappa^2/3)}\Big].$$

\end{lemma}
Thanks to Lemma \ref{REn}, the empirical Restricted Eigenvalue condition will be fulfilled on an event of large probability, on which we establish a fast non-asymptotic oracle inequality.


\begin{theorem}
\label{thm:fast-oracle}
Consider Model (\ref{eq:CoxGen}) with known $\alpha_{0}$ and for $x>0$, let $\omega_{j}$ be defined by (\ref{eq:weights}) and $\bs{\hat\beta_L^\mu}$ be defined by (\ref{eq:lassoGamma}). Let $A_{\varepsilon,\nu}>0$ be a numerical positive constant only depending on $\varepsilon$ and $\nu$, $\zeta>0$ and $s\in\{1,...,M\}$ be fixed. Let Assumptions \ref{ass:dico1}, \ref{ass:voisinage} and \ref{ass:RE} be satisfied with $a_{0}= (3+4/\zeta)$ and let $\bs\kappa=(1/\sqrt{2A_0})\bs\kappa_0(s,a_{0})$.
Then, with a probability larger than $1-~A_{\varepsilon,\nu}\e^{-x}-\pi_n$, the following inequality holds
\begin{equation}
\label{eq:fast-oracle}
\widetilde{K}_{n}(\lambda_{0},\lambda_{\hat{\beta}_{L}^{\mu}})\leq(1+\zeta)\underset{|J(\bs{\beta})|\leq s}{\underset{\bs{\beta}\in\Gamma(\mu)}{\inf}}\left\{\widetilde{K}_{n}(\lambda_{0},\lambda_{\beta})+C(\zeta,\mu)\dfrac{|J(\bs{\beta})|}{\bs\kappa^2}(\underset{1\leq j\leq M}{\max}\omega_{j})^2\right\},
\end{equation}
where $C(\zeta,\mu)>0$ is a constant depending on $\zeta$ and $\mu$.
\end{theorem}
This result allows to compare the prediction error of the estimator and the best sparse approximation of the regression function by an oracle that knows the truth, but is constrained by sparsity. The Lasso estimator approaches the best approximation in the dictionary with a fast error term of order $\log M/n$. 

Thanks to Proposition \ref{prop:comparaison}, which states a connection between the empirical Kullback divergence (\ref{eq:Kullback}) and the weighted empirical quadratic norm (\ref{eq:norm}), we deduce from Theorem \ref{thm:fast-oracle} a non-asymptotic oracle inequality in weighted empirical quadratic norm.
 \begin{corollary}
\label{cor:corollary}
Under the assumptions of Theorem \ref{thm:fast-oracle}, with a probability larger than $1-A_{\varepsilon,\nu}\e^{-x}-\pi_n$,
\begin{equation*}
||\log\lambda_{\hat\beta_{L}^{\mu}}-\log\lambda_{0}||_{n,\Lambda}^2\leq(1+\zeta)\underset{|J(\beta)|\leq s}{\underset{\beta\in\Gamma(\mu)}{\inf}}\left\{||\log\lambda_{\beta}-\log\lambda_{0}||_{n,\Lambda}^2+\tilde c(\zeta,\mu)\dfrac{|J(\bs{\beta})|}{\bs\kappa^2}(\underset{1\leq j\leq M}{\max}\omega_{j})^2\right\},
\end{equation*}
where $\tilde c(\zeta,\mu)$ is a positive constant depending on $\zeta$ and $\mu$.
\end{corollary}
Note that for $\alpha_{0}$ supposed to be known, this oracle inequality is also equivalent to
\begin{equation*}
||f_{\hat\beta_{L}^{\mu}}-f_{0}||_{n,\Lambda}^2\leq(1+\zeta)\underset{|J(\beta)|\leq s}{\underset{\beta\in\Gamma(\mu)}{\inf}}\left\{||f_{\beta}-f_{0}||_{n,\Lambda}^2+\tilde c(\zeta,\mu)\dfrac{|J(\bs{\beta})|}{\bs\kappa^2}(\underset{1\leq j\leq M}{\max}\omega_{j})^2\right\}.
\end{equation*}

\subsection{Particular case : variable selection in the Cox model}

We now consider the case of variable selection in the Cox model (\ref{eq:CoxGen}) with $f_{0}(Z_{i})=\bs{\beta_{0}^TZ_{i}}$.
In this case, $M=p$ and the functions of the dictionary are such that for $i=1,...,n$ and $j=1,...,p$
$$f_{j}(\bs{Z_{i}})=Z_{i,j} \text{ and } f_{\beta}(\bs{Z_{i}})=\sum_{j=1}^{p}\beta_{j}Z_{i,j}=\bs{\beta^TZ_{i}}.$$

Let $\bs{X}=(Z_{i,j})_{\underset{1\leq j\leq p}{1\leq i\leq n}}$ be the design matrix and for $\bs{\hat\beta_{L}}$ defined by (\ref{eq:lassosimple}), let

$$\bs{\Delta_{0}}=\bs{D(\hat\beta_{L}}-\bs{\beta_{0})}, \bs{D}=(\diag(\omega_j))_{1\leq j\leq M}, J_{0}=J(\bs{\beta_{0}}) \text{ and } |J_{0}|=\Card\{J_{0}\}.$$

We now state non-asymptotic inequalities for prediction on $\bs{X\beta_{0}}$ and for estimation on $\bs{\beta_{0}}$. 
In this subsection, we don't need to work locally on the set $\Gamma(\mu)$ to obtain Proposition \ref{relationsel} and instead of considering Assumption (\ref{ass:voisinage}), we only have to introduce the following assumption to connect the empirical Kullback divergence and the weighted empirical quadratic norm :
\begin{align}
&\hyp{covborne}  \text{Let } R \text{ be a positive constant, such that } \underset{i\in\{1,...,n\}}{\max}||Z_i||_{2}\leq R. \notag
\end{align}

We consider the Lasso estimator defined with the regularization parameter $\Gamma_1>0$ :
$$\bs{\hat\beta_L}=\underset{\bs\beta\in\mathbb{R}^p}{\arg\min}\{C_n(\lambda_{\bs\beta})+\Gamma_1\displaystyle\sum_{j=1}^p\omega_j|\beta_j|\},$$

\begin{theorem}
\label{thm:selection}
Consider Model (\ref{eq:Cox}) with known $\alpha_{0}$. For $x>0$, let $\omega_{j}$ be defined by (\ref{eq:weights}) and denote $\bs{\kappa'}=(1/\sqrt{2A_0})\bs{\kappa_0}(s,3).$
Let $A_{\varepsilon,\nu}$ be some numerical positive constant depending on $\varepsilon$ and $\nu$. Under Assumptions \ref{ass:dico1}, \ref{covborne} and \ref{ass:RE} with $a_0=3$, for all $\Gamma_1$ such that
$$\Gamma_1\leq \dfrac{1}{48Rs}\dfrac{\underset{1\leq j\leq M}{\min}\omega_j^2}{\underset{1\leq j\leq M}{\max}\omega_j^2}\dfrac{\bs\kappa'^{2}}{\underset{1\leq j\leq M}{\max}\omega_j},$$ 
with a probability larger than $1-A_{\varepsilon,\nu}\e^{-\Gamma_1x}-\pi_n$, then
 \begin{equation}
\label{eq:prediction}
||\bs{X}(\bs{\hat\beta_{L}}-\bs{\beta_{0}})||_{n,\Lambda}^2\leq \dfrac{4}{\xi^2} \dfrac{|J_{0}|}{\bs{\kappa'}^2}\Gamma_1^2(\underset{1\leq j\leq p}{\max}\omega_{j})^2
\end{equation}
and 
\begin{equation}
\label{eq:selection}
||\bs{\hat\beta_{L}}-\bs{\beta_{0}}||_{1}\leq 8\dfrac{\underset{1\leq j\leq p}{\max}\omega_j}{\underset{1\leq j\leq p}{\min}\omega_j}\dfrac{|J_{0}|}{\xi\bs{\kappa'}^2}\Gamma_1\underset{1\leq j\leq p}{\max}\omega_{j}.
\end{equation}
\end{theorem}

This theorem gives non-asymptotic upper bounds for two types of loss functions. Inequality (\ref{eq:prediction}) gives a non-asymptotic bound on prediction loss with a rate of convergence in ${\log M/n}$, while Inequality (\ref{eq:selection}) states a bound on $\bs{\hat\beta_{L}}-\bs{\beta_{0}}$. 


\section{Oracle inequalities for general intensity}

In the previous section, we have assumed $\alpha_0$ known and have obtained results on the relative risk. Now, we consider a general intensity $\lambda_{0}$ that does not rely on an underlying model. Oracle inequalities are established under different assumptions with slow and fast rates of convergence.

\subsection{A slow oracle inequality}
\label{subsec:slow-inconnu}
The slow oracle inequality for a general intensity is obtained under light assumptions that concern only the construction of the two dictionaries $\mathbb{F}_{M}$ and $\mathbb{G}_{N}$.
\begin{theorem}
\label{thm:inconnu}
For $x>0$ and $y>0$, let $\omega_{j}$ and $\delta_k$ be defined by (\ref{eq:weights}) and $(\bs{\hat\beta_L},\bs{\hat\gamma_L})$ be defined in Estimation procedure \ref{EstimationProc}.
Let $A_{\varepsilon,\nu}$ and $B_{\tilde\varepsilon,\tilde\nu}>0$ be two positive numerical constants depending on $\varepsilon, \nu \text{ and } \tilde\varepsilon, \tilde\nu $ respectively and Assumptions \ref{ass:dico1}, \ref{ass:dico2} be satisfied. Then, with probability larger than $1-A_{\varepsilon,\nu}\e^{-x}-B_{\tilde\varepsilon,\tilde\nu}\e^{-y}$
\begin{align}
\label{eq:slow-inconnu}
\widetilde{K}_n(\lambda_{0},\lambda_{\hat\beta_{L},\hat\gamma_{L}})\leq\underset{(\bs\beta,\bs\gamma)\in\mathbb{R}^M\times\mathbb{R}^N}{\inf}\{\widetilde{K}_n(\lambda_{0},\lambda_{\beta,\gamma})+2\pen(\bs\beta)+2\pen(\bs\gamma)\}.
\end{align}
\end{theorem}

We have chosen to estimate the complete intensity, which involves two different parts~ : the first part is the baseline function $\alpha_{\bs\gamma} : \mathbb{R}\rightarrow\mathbb{R}$ and the second part is the function of the covariates $f_{\bs\beta} : \mathbb{R}^p\rightarrow\mathbb{R}$. The double $\ell_1$-penalization considered here is tuned to concurrently estimate the function $f_0$ depending on high-dimensional covariates and the non-parametric function $\alpha_0$. Examples of Lasso algorithms for the estimation of non-parametric density or intensity may be found in Bertin et al. \cite{BlPR} and Hansen et al. \cite{Hansen} respectively. As $f_0$ and $\alpha_0$ are estimated at once, the resulting rate of convergence is the sum of the two expected rates in both situations considered separately ($\sim\sqrt{\log M/n}+\sqrt{\log N/n}$). Nevertheless, from Bertin et al. \cite{BlPR}, we expect that a choice of $N$ of order $n$ would suitably estimate $\alpha_0$. As a consequence, in a very high-dimensional setting the leading error term in (\ref{eq:slow-inconnu}) would be of order $\sqrt{\log M/n}$, which again is the classical slow rate of convergence in a regression setting.

\subsection{A fast oracle inequality}
\label{subsec:fast-inconnu}
We are now interested in obtaining the fast non-asymptotic oracle inequality and as usual, we need to introduce further notations and assumptions.
In this subsection, we shall again work locally for $\rho>0$ on the set $\widetilde\Gamma_{M,N}(\rho)=\{(\bs\beta,\bs\gamma)\in\mathbb{R}^M\times\mathbb{R}^N : ||\log\lambda_{\bs\beta,\bs\gamma}-~\log\lambda_0||_{n,\infty}\leq \rho\}$, simply denoted $\widetilde\Gamma(\rho)$ and we consider the following assumption :

\begin{align}
&\hyp{ass:voisinage1}  \text{There exists } \rho>0, \text{ such that } \widetilde\Gamma(\rho) \text{ contains a non-empty open set of } \mathbb{R}^M\times\mathbb{R}^N. \notag
\end{align}
On $\widetilde\Gamma(\rho)$, we define the weighted Lasso estimator as   
$$(\bs{\hat\beta_L^{\rho}},\bs{\hat\gamma_L^{\rho}})=\underset{(\bs\beta,\bs\gamma)\in\widetilde\Gamma(\rho)}{\arg\min}\{C_n(\lambda_{\bs\beta,\bs\gamma})+\pen(\bs\beta)+\pen(\bs\gamma)\}.$$

Let us give the additional notations. Set $\bs{\tilde\Delta}$ be
$$\bs{\tilde\Delta}=\bs{\tilde D}\begin{pmatrix} \bs{\hat\beta_{L}}-\bs\beta\\
                                                                  \bs{\hat\gamma_{L}}-\bs\gamma
\end{pmatrix} \in \mathbb{R}^{M+N} \text{ with } (\bs\beta,\bs\gamma)\in\widetilde\Gamma(\rho) \text{ and }\tilde{\bs{D}}=\diag(\omega_{1},...,\omega_{M},\delta_{1},...,\delta_{N}).$$
Let $\bs{\Large{1}}_{n\times N}$ be the matrix $n\times N$ with all coefficients equal to one,

%

%

$$\bs{\tilde X}(t)=\begin{bmatrix} (f_{j}(\bs{Z_{i}}))_{\underset{1\leq j\leq M}{1\leq i\leq n}} & \bs{\Large{1}}_{n\times N}(\diag(\theta_{k}(t)))_{1\leq k\leq N}\end{bmatrix}=\left[\begin{array}{c|ccc} & \theta_{1}(t) & \dots &\theta_{N}(t)\\
                                                        \bs X   & \vdots& &\vdots\\
                                                        \text{}  &  \theta_{1}(t)&\dots  & \theta_{N}(t) 

\end{array}\right]\in\mathbb{R}^{n\times(M+N)}$$

%

and
$$\bs{\tilde{G}_{n}}=\dfrac{1}{n}\displaystyle{\inttau\bs{\tilde{X}}(t)^T\bs{\tilde{C}}(t)\bs{\tilde{X}}(t)\diff t}\text{ with }\bs{\tilde C}(t)=(\diag(\lambda_{0}(t,\bs{Z_{i}})Y_{i}(t)))_{1\leq i\leq n}, \forall t\geq 0.$$

Let also $J(\bs{\beta})$ and $J(\bs\gamma)$ be the sparsity sets  of vectors $(\bs{\beta},\bs\gamma)\in\widetilde\Gamma(\rho)$ respectively defined by 
$$J(\bs{\beta})=\{j\in\{1,...,M\} : \beta_{j}\neq0\} \text{ and }J(\bs{\gamma})=\{k\in\{1,...,N\} : \gamma_{k}\neq0\},$$
and the sparsity indexes are then given by
$$|J(\bs{\beta})|=\sumip\mathds{1}_{\{\beta_{j}\neq 0\}}=\Card\{J(\bs\beta)\} \text{ and } |J(\bs{\gamma})|=\displaystyle\sum_{k=1}^{N}\mathds{1}_{\{\gamma_{k}\neq 0\}}=\Card\{J(\bs\gamma)\}.$$
To obtain the fast non-asymptotic oracle inequality, we consider the Restricted Eigenvalue condition applied to the matrix $\mathbb{E}(\bs{\tilde{G}_n})$.
\begin{align}
&  \text{For some integer } s\in\{1,...,M+N\} \text{ and a constant } r_{0}>0, \text{ we assume that } \bs{\tilde{G}_{n}} \text{ satisfies :} \notag\\
&\hspace{2cm}\hypREinconnu{ass:RE1} 0<\bs{\tilde\kappa_0}(s,r_{0})=\underset{|J|\leq s}{\underset{J\subset\{1,...,M+N\},}{\min}}\underset{||\bs{b}_{J^c}||_{1}\leq r_{0}||\bs{b}_{J}||_{1}}{\underset{\bs{b}\in\mathbb{R}^{M+N}\backslash\{0\},}{\min}}\dfrac{(\bs{b^T\mathbb{E}(\tilde{G}_{n})b})^{1/2}}{||\bs{b}_{J}||_{2}}. \notag 
\end{align}

The condition on the matrix $\mathbb{E}(\bs{\tilde{G}_{n}})$ is rather strong because the block matrix involves both functions of the covariates of $\mathbb{F}_{M}$ and functions of time which belong to $\mathbb{G}_{N}$. This is the price to pay for an oracle inequality on the full intensity.  If we had instead considered two restricted eigenvalue assumptions on each block, we would have established an oracle inequality on the sum of the two unknown parameters $\alpha_{0}$ and $f_{0}$ and not on $\lambda_{0}$. As in Lemma \ref{REn}, we can show that under Assumption \ref{ass:RE1}, we have an empirical Restricted Eigenvalue condition on the matrix $\bs{\tilde{G}_n}$.

\begin{lemma}
\label{tildeREn}
Let $L$ defined as in Lemma \ref{REn}. Under Assumptions \ref{A_0} and \ref{ass:RE1}, we have 
\begin{equation}
\label{eq:tildeREn}
0<\bs{\tilde\kappa}=\underset{|J|\leq s}{\underset{J\subset\{1,...,M\},}{\min}}\underset{||\bs{b}_{J^c}||_{1}\leq r_{0}||\bs{b}_{J}||_{1}}{\underset{\bs{b}\in\mathbb{R}^M\backslash\{0\},}{\min}}\dfrac{(\bs{b^T\tilde{G}_nb})^{1/2}}{||\bs{b}_{J}||_{2}} \mbox{ and } \bs{\tilde\kappa}=(1/\sqrt{2A_0})\bs{\tilde{\kappa}_0}(s,r_0),
\end{equation}
with probability larger than $1-\tilde\pi_n$, where
$$\tilde\pi_n=2M^2\exp\Big[-\dfrac{n\bs{\tilde\kappa}^4}{2L^2(1+r_0)^2s(L^2(1+r_0)^2s+\bs{\tilde\kappa}^2/3)}\Big].$$

\end{lemma}


\begin{theorem}
\label{thm:fast-inconnu}
For $x>0$ and $y>0$, let $\omega_{j}$ and $\delta_{k}$ be defined by (\ref{eq:weights}). Let $A_{\varepsilon,\nu}>0 \text{ and } B_{\tilde\varepsilon,\tilde\nu}>0$ be two numerical positive constants depending on $\varepsilon, \nu \text{ and } \tilde\varepsilon, \tilde\nu$ respectively, $\zeta>0$ and $s\in\{1,...,M+N\}$ be fixed. Let Assumptions \ref{ass:dico1}, \ref{ass:dico2}, \ref{ass:voisinage1} and \ref{ass:RE1} be satisfied with
$$r_{0}= \left(3+8\max\left(\sqrt{|J(\bs\beta)|},\sqrt{|J(\bs\gamma)|}\right)/\zeta\right),$$ 
and let $\bs{\tilde\kappa}=(1/\sqrt{2A_0})\bs{\tilde\kappa_0}(s,r_{0})$.
Then, with probability larger than $1-A_{\varepsilon,\nu}\e^{-x}-B_{\tilde\varepsilon,\tilde\nu}\e^{-y}-\tilde\pi_n$

\begin{align}
\label{eq:fast-inconnu}
\widetilde{K}_{n}(\lambda_{0},\lambda_{\hat\beta_{L}^{\rho},\hat\gamma_{L}^{\rho}})\leq (1+\zeta)\underset{\max(|J(\bs\beta)|,|J(\bs\gamma)|)\leq s}{\underset{(\bs\beta,\bs\gamma)\in\widetilde\Gamma(\rho)}{\inf}}\Big\{\widetilde{K}_{n}(\lambda_{0},\lambda_{\beta,\gamma})+\widetilde{C}(\zeta,\rho)\dfrac{\max(|J(\bs\beta)|,|J(\bs\gamma)|)}{\bs{\tilde\kappa}^2}\underset{1\leq k\leq N}{\underset{1\leq j\leq M}{\max}}\{\omega_{j}^2,\delta_{k}^2\}\Big\},
\end{align}
and
\begin{align}
\label{eq:fast-inconnu-emp}
\nonumber
||\log&\lambda_{0}-\log\lambda_{\hat\beta_{L}^{\rho},\hat\gamma_{L}^{\rho}}||_{n,\Lambda}^2\\
&\leq (1+\zeta)\underset{\max(|J(\bs\beta)|,|J(\bs\gamma)|)\leq s}{\underset{(\bs\beta,\bs\gamma)\in\widetilde\Gamma(\rho)}{\inf}}\Big\{||\log\lambda_{0}-\log\lambda_{\beta,\gamma}||_{n,\Lambda}^2+\widetilde{C}'(\zeta,\rho)\dfrac{\max(|J(\bs\beta)|,|J(\bs\gamma)|)}{\bs{\tilde\kappa}^2}\underset{1\leq k\leq N}{\underset{1\leq j\leq M}{\max}}\{\omega_{j}^2,\delta_{k}^2\}\Big\},
\end{align}
where $\widetilde{C}(\zeta,\rho)>0$ and $\widetilde{C}'(\zeta,\rho)>0$ are constants depending only on $\zeta$ and $\rho$.
\end{theorem}
We obtain a non-asymptotic fast oracle inequality in prediction. Indeed, the rate of convergence of this oracle inequality is of order
$$\Big(\underset{1\leq k\leq N}{\underset{1\leq j\leq M}{\max}}\{\omega_{j},\delta_{k}\}\Big)^2\approx \max\Big\{\dfrac{\log M}{n},\dfrac{\log N}{n}\Big\},$$
namely, if we choose $\mathbb{G}_{N}$ of size $n$, the rate of convergence of this oracle inequality is then of order $\log M/n$ (see Subsection \ref{subsec:slow-inconnu} for more details). 
While Estimation procedure \ref{EstimationProc} allows to derive a prediction for the survival time through the conditional intensity, Theorem \ref{thm:fast-inconnu} measures the accuracy of this prediction. In that sense, the clinical problem of establishing a prognosis has been addressed at this point. To our best knowledge, this oracle inequality is the first non-asymptotic oracle inequality in prediction for the whole intensity with a fast rate of convergence of order $\log M/n$. 


For the part depending on the covariates, recent results establish non-asymptotic oracle inequalities for the Lasso estimator of $f_0$ in the usual Cox model (see Bradic and Song \cite{bradic2012} and Kong and Nan \cite{Kong}). We cannot compare our results to theirs, since we estimate the whole intensity with the total empirical log-likelihood whereas both of them consider the partial log-likelihood. 


The remaining part of the paper is devoted to the technical results and proofs

\section{An empirical Bernstein's inequality}

The main ingredient of Theorems \ref{thm:slow-oracle}, \ref{thm:fast-oracle}, \ref{thm:inconnu} and \ref{thm:fast-inconnu} are Bernstein's concentration inequalities that we present in this section. To clarify the relation between the stated oracle inequalities and the Bernstein's inequality, we sketch here the proof of Theorem \ref{thm:inconnu}.
Using the Doob-Meyer decomposition $N_{i}=M_{i}+\Lambda_{i}$, we can easily show that for all $\bs\beta\in\mathbb{R}^M$ and for all $\bs\gamma\in\mathbb{R}^N$
\begin{equation}
\label{eq:difference-risk}
C_{n}(\lambda_{\hat\beta_{L},\hat\gamma_{L}})-C_{n}(\lambda_{\beta,\gamma})=\widetilde{K}_{n}(\lambda_{0},\lambda_{\hat\beta_{L},\hat\gamma_{L}})-\widetilde{K}_{n}(\lambda_{0},\lambda_{\beta,\gamma})+(\bs{\hat\gamma_{L}}-\bs\gamma)^T\bs{\nu_{n,\tau}}+(\bs{\hat\beta_{L}}-\bs\beta)^T\bs{\eta_{n,\tau}},
\end{equation}
where 

\begin{align}
\label{eq:proctocontrol}
\bs{\eta_{n,\tau}}=\dfrac{1}{n}\sumin\inttau\bs{\vec{f}}(\bs{Z_{i}})\diff M_{i}(t) \text{ andÊ} \bs{\nu_{n,\tau}}=\dfrac{1}{n}\sumin\inttau\bs{\vec{\theta}}(t)\diff M_{i}(t), 
\end{align}
with $\bs{\vec{f}}=(f_{1},...,f_{M})^T$ and $\bs{\vec{\theta}}=(\theta_{1},...,\theta_{N})^T$.
By definition of the Lasso estimator, we have for all $(\bs\beta,\bs\gamma)$ in $\mathbb{R}^M\times\mathbb{R}^N$ 
$$C_{n}(\lambda_{\hat\beta_{L},\hat\gamma_{L}})+\pen(\bs{\hat\beta_L})+\pen(\bs{\hat\gamma_L})\leq C_{n}(\lambda_{\beta,\gamma})+\pen(\bs\beta)+\pen(\bs\gamma),$$
and we finally obtain 
$$\widetilde{K}_{n}(\lambda_{0},\lambda_{\hat\beta_{L},\hat\gamma_{L}})\leq \widetilde{K}_{n}(\lambda_{0},\lambda_{\beta,\gamma})+(\bs{\hat\gamma_{L}}-\bs\gamma)^T\bs{\nu_{n,\tau}}+(\bs{\hat\beta_{L}}-\bs\beta)^T\bs{\eta_{n,\tau}}+\pen(\bs\beta)-\pen(\bs{\hat\beta_L})+\pen(\bs\gamma)-\pen(\bs{\hat\gamma_L}).$$
 Consequently, $\widetilde{K}_{n}(\lambda_{0},\lambda_{\hat\beta_{L},\hat\gamma_{L}})$ is bounded by
\begin{align*}
\widetilde{K}_{n}(\lambda_{0},\lambda_{\beta,\gamma})+\displaystyle\sum_{j=1}^{M}({\hat\beta_{L,j}}-\beta_j)\eta_{n,\tau}(f_j)+\displaystyle\sum_{j=1}^M\omega_j(|\beta_j|-|\hat\beta_{L,j}|)+\displaystyle\sum_{k=1}^{N}({\hat\gamma_{L,k}}-\gamma_k)^T{\nu_{n,\tau}(\theta_k)}+\displaystyle\sum_{k=1}^N\delta_k(|\gamma_k|-|\hat\gamma_{L,k}|),
\end{align*}
with
\begin{align*}
\eta_{n,t}(f_{j})=\dfrac{1}{n}\displaystyle{\sum_{i=1}^{n}\int_{0}^{t}f_{j}(\bs{Z_{i}})\diff M_{i}(s)} \text{ and } \nu_{n,t}(\theta_{k})=\dfrac{1}{n}\sumin\intt\theta_{k}(s)\diff M_{i}(s).
\end{align*}
We will control $\eta_{n,t}(f_{j})$ and $\nu_{n,t}(\theta_{k})$ respectively by $\omega_j$ and $\delta_k$. More precisely, the weights $\omega_j$ (respectively $\delta_k$) will be chosen 
such that $|\eta_{n,t}(f_{j})|\leq \omega_j$ (respectively $|\nu_{n,t}(\theta_{k})|\leq \delta_k$) and $\mathbb{P}(|\eta_{n,t}(f_{j})|> \omega_j)$ (respectively $\mathbb{P}(|\nu_{n,t}(\theta_{k})|> \delta_k$) large. As $\eta_{n,t}(f_j)$ and $\nu_{n,t}(\theta_k)$ involve martingales, we could directly apply classical Bernstein's inequalities for martingales with $x>0$ and $y>0$ 
$$\mathbb{P}\Big[\eta_{n,t}(f_{j})\geq \sqrt{\dfrac{2V_{n,t}(f_{j}) x}{n}}+\dfrac{x}{3n}\Big]\leq \e^{-x} \text{ and }\mathbb{P}\Big[\nu_{n,t}(\theta_{k})\geq \sqrt{\dfrac{2R_{n,t}(\theta_{k}) y}{n}}+\dfrac{y}{3n}\Big]\leq \e^{-y} ,$$
where the predictable variations $V_{n,t}(f_{j})$ and $R_{n,t}(\theta_{k})$ of $\eta_{n,t}(f_{j})$ and $\nu_{n,t}(\theta_{k})$ are respectively defined by
 \begin{align*}
V_{n,t}(f_{j})&=n<\eta_{n}(f_{j})>_{t}=\dfrac{1}{n}\displaystyle{\sum_{i=1}^{n}\int_{0}^{t}(f_{j}(\bs{Z_{i}}))^2\lambda_{0}(t,\bs{Z_{i}})Y_{i}(s)\diff s},\\
R_{n,t}(\theta_{k})&=n<\nu_{n}(\theta_{k})>_{t}=\dfrac{1}{n}\sumin\intt(\theta_{k}(t))^2\lambda_{0}(t,\bs{Z_{i}})Y_{i}(s)\diff s,
\end{align*}
see e.g. van de Geer \cite{SVG}.
Applying these inequalities, the weights of Algorithm \ref{EstimationProc} would have the forms $\omega_j=\sqrt{2V_{n,t}(f_{j}) x/n}+x/3n$ and $\delta_k=\sqrt{2R_{n,t}(\theta_{k}) y/n}+y/3n$.
As $V_{n,t}(f_{j})$ and $R_{n,t}(\theta_{k})$ both depend on $\lambda_0$, this would not result a statistical procedure. We propose to replace in the Bernstein's inequality the predictable variations by the optional variations of the processes $\eta_{n,t}(f_{j})$ and $\nu_{n,t}(\theta_{k})$ defined by
\begin{align*}
\hat{V}_{n,t}(f_{j})=n[\eta_{n}(f_{j})]_{t}=\dfrac{1}{n}\displaystyle{\sum_{i=1}^{n}\int_{0}^{t}(f_{j}(\bs{Z_{i}}))^2\diff N_{i}(s)} \text{ and } \hat{R}_{n,t}(\theta_{k})=n[\nu_{n}(\theta_{k})]_{t}=\dfrac{1}{n}\sumin\intt(\theta_{k}(t))^2\diff N_{i}(s).
\end{align*}
This ensures that the weights $\omega_j$ and $\delta_k$ will depends on $\hat{V}_{n,t}(f_{j})$ and $\hat{R}_{n,t}(\theta_{k})$ respectively. Equivalent strategies in different models have been considered in Ga\"iffas and Guilloux \cite{Agathe} or Hansen et al. \cite{Hansen}. The following theorem states the resulting Bernstein's inequalities.

\begin{theorem}
\label{thm:empirical-bernstein}
Let Assumption \ref{A_0} be satisfied. For any numerical constant $\varepsilon>0$, $\tilde\varepsilon>0$, $c=\sqrt{2(1+\varepsilon)}$ and $\tilde c=\sqrt{2(1+\tilde\varepsilon)}$,  the following holds for any $x>0$, $y>0$ :
\begin{align}
\label{eq:empirical-bernstein}
\mathbb{P}\Big[|\eta_{n,t}(f_{j})|\geq c\sqrt{\dfrac{\hat W_n^{\nu}(f_j)x}{n}}+\dfrac{x}{3n}||f_j||_{n,\infty}\Big]\leq \Big(\dfrac{2}{\log(1+\varepsilon)}\log\Big(2+\dfrac{A_0(\nu/n+\Phi(\nu/n))}{x/n}\Big)+1\Big)\e^{-x},\\
\label{eq:empirical-bernstein1}
\mathbb{P}\Big[|\nu_{n,t}(\theta_{k})|\geq \tilde c\sqrt{\dfrac{\hat T_n^{\tilde\nu}(\theta_{k})y}{n}}+\dfrac{y}{3n}||\theta_{k}||_{\infty}\Big]\leq \Big(\dfrac{2}{\log(1+\tilde\varepsilon)}\log\Big(2+\dfrac{A_0(\tilde\nu/n+\Phi(\tilde\nu/n))}{y/n}\Big)+1\Big)\e^{-y},
\end{align}
where
\begin{align}
\label{foncopt}
W_n^{\nu}(f_j)&=\dfrac{\nu/n}{\nu/n-\Phi(\nu/n)}\hat V_n(f_j)+\dfrac{x/n}{\nu/n-\Phi(\nu/n)}||f_j||_{n,\infty}^2,\\
T_n^{\tilde\nu}(\theta_k)&=\dfrac{\tilde\nu/n}{\tilde\nu/n-\Phi(\tilde\nu/n)}\hat R_n(\theta_k)+\dfrac{y/n}{\tilde\nu/n-\Phi(\tilde\nu/n)}||\theta_k||_{\infty}^2, 
\end{align}
for real numbers $(\nu,\tilde\nu)\in(0,3)^2$ such that $\nu>\Phi(\nu)$ and $\tilde\nu>\Phi(\tilde\nu)$, where $\Phi(u)=\exp(u)-u-1$.
\end{theorem}
We deduce the weights $\omega_j$ and $\delta_k$ defined in (\ref{eq:weights}), from Theorem \ref{thm:empirical-bernstein}.
These empirical Bernstein's inequalities hold true for martingales with jumps, when the predictable variation is not observable. 

\begin{remark}
Theorem \ref{thm:empirical-bernstein} is closed to Theorem 3 in Hansen et al. \cite{Hansen}, although in our version the event bounding $\hat W_n^{\nu}(f_j)$ and $\hat T_{n}^{\tilde\nu}(\theta_k)$ has been removed from the probability (see the proof of Theorem  \ref{thm:empirical-bernstein}).

Other weights can also be obtained from empirical Bernstein's inequalities that are closer to those obtained by Ga\"iffas and Guilloux \cite{Agathe} in Theorem 3. We refer to an other version of the paper (see \cite{Lemler12}), in which these weights appear. Their forms are less simple than those defined in (\ref{eq:weights}), but they do not depend on tuning parameters $\nu$ and $\tilde\nu$ to determine for the applications. An interesting perspective would be to determine which one of those two forms of weights gives the best results in the applications.
\end{remark}

\subsection*{Acknowledgements}
All my thanks go to my two Phd Thesis supervisors Agathe Guilloux and Marie-Luce Taupin for their help, their availability and their advices. I also thank Marius Kwemou for helpful discussions.

\section{Proofs}

\subsection{Proof of Proposition \ref{prop:Kullback}}

Following the proof of Theorem 1 in Senoussi \cite{Senoussi}, we rewrite the empirical Kullback divergence (\ref{eq:Kullback}) as 
\begin{align*}
\widetilde{K}_{n}(\lambda_{0},\lambda)&=\dfrac{1}{n}\displaystyle{\sum_{i=1}^{n}\inttau\Big[\log\lambda_{0}(t,\bs{Z_{i}})-\log\lambda(t,\bs{Z_{i}})-\Big(1-\dfrac{\lambda(t,\bs{Z_{i}})}{\lambda_{0}(t,\bs{Z_{i}})}\Big)\Big]\lambda_{0}(t,\bs{Z_{i}})Y_{i}(t)\diff t}\\
&=\dfrac{1}{n}\displaystyle{\sum_{i=1}^{n}\inttau\Big[\exp\Big({\log{\frac{\lambda(t,\bs{Z_{i}})}{\lambda_{0}(t,\bs{Z_{i}})}}}\Big)-\log\dfrac{\lambda(t,\bs{Z_{i}})}{\lambda_{0}(t,\bs{Z_{i}})}-1\Big]\lambda_{0}(t,\bs{Z_{i}})Y_{i}(t)\diff t}.
\end{align*}
Since the map $t\rightarrow \e^t-t-1$ is a positive function on $\mathbb{R}$, we deduce that except for $\lambda=\lambda_{0},$
$$\exp\Big({\log{\frac{\lambda(t,\bs{Z_{i}})}{\lambda_{0}(t,\bs{Z_{i}})}}}\Big)-\log\dfrac{\lambda(t,\bs{Z_{i}})}{\lambda_{0}(t,\bs{Z_{i}})}-1>0.$$ 
Thus $\widetilde{K}_{n}(\lambda_{0},\lambda)$ is positive and vanishes only if  $(\log\lambda_{0}-\log\lambda)(t,\bs{Z_{i}})=0$ almost surely, namely if $\lambda_{0}=\lambda$ almost surely. 
\qed

\subsection{Proof of Proposition \ref{thm:slow-oracle}}

According to the definition (\ref{eq:lassoGamma}) of $\bs{\hat{\beta}_{L}}$, for all $\beta$ in $\mathbb{R}^M$, we have 
\begin{equation*}
C_{n}(\lambda_{\hat{\beta}_{L}})+\pen(\bs{\hat{\beta}_{L}})\leq C_{n}(\lambda_{\beta})+\pen(\bs{\beta}).
\end{equation*}
Here $\alpha_0$ is assumed to be known. Hence applying (\ref{eq:difference-risk}), we obtain
\begin{align}
\label{eq:B1}
\widetilde{K}_{n}(\lambda_{0},\lambda_{\hat\beta_{L}})\leq\widetilde{K}_{n}(\lambda_{0},\lambda_{\beta})+(\bs{\hat{\beta}_{L}}-\bs{\beta})^{T}\bs{\eta_{n,\tau}}+\pen(\bs{\beta})-\pen(\bs{\hat{\beta}_{L}}).
\end{align} 
It remains to control the term $(\bs{\hat{\beta}_{L}}-\bs{\beta})^{T}\bs{\eta_{n,\tau}}$.
For $\omega_j$ defined in (\ref{eq:weights}), set
\begin{align*}
\mathcal{A}=\underset{j=1}{\overset{M}{\bigcap}}\left\{|\eta_{n,\tau}(f_{j})|\leq\dfrac{\omega_{j}}{2}\right\}. 
\end{align*}
On $\mathcal{A}$, we have 
$$|(\bs{\hat{\beta}_{L}}-\bs\beta)^T\bs{\eta_{n,\tau}}|\leq\sumip\dfrac{\omega_{j}}{2}|(\hat{\beta}_{L}-\beta)_{j}|\leq\sumip\omega_{j}|(\hat{\beta}_{L}-\beta)_{j}|.$$
The result (\ref{eq:slow-oracle}) follows since $\pen(\bs\beta)=\displaystyle\sumip\omega_{j}|\beta_{j}|$.  
It remains to bound up $\mathbb{P}(\mathcal{A}^c)$. By applying Theorem \ref{thm:empirical-bernstein}
\begin{align*}
\mathbb{P}(\mathcal{A}^c)&\leq\sumip\mathbb{P}\Big(|\eta_{n,\tau}(f_{j})|>\dfrac{\omega_{j}}{2}\Big)\leq A_{\varepsilon,\nu}\e^{-x},
\end{align*}
with
$$A_{\varepsilon,\nu}=\dfrac{2}{\log(1+\varepsilon)}\log\Big(2+\dfrac{A_0(\nu/n+\Phi(\nu/n))}{x/n}\Big)+1.$$
We conclude that $\mathbb{P}(\mathcal{A})\geq1-A_{\varepsilon,\nu}\e^{-x},$
which ends up the proof of Theorem \ref{thm:slow-oracle}.
\qed 

\subsection{Proof of Lemma \ref{REn}}
We show with high probability, that under \ref{ass:RE}, for all $J\subset\{1,...,M\}$ such that $|J|\leq s$ and for all $b\in \mathbb{R}^M\backslash\{0\}$ such that $||b_{J^c}||_1\leq a_0||b_J||_1$,
$$\dfrac{\bs{b^TG_nb}}{||\bs{b}_J||_2^2}>\bs\kappa^2,\quad \mbox{with} \quad \bs\kappa=(1/\sqrt{2A_0})\bs{\kappa_0}(s,a_0) \quad \mbox{and} \quad A_0 \mbox{ defined in Assumption \ref{A_0}}.$$
Let consider the set $\Omega_{G_n}=\{|(\bs{G_n}-\mathbb{E}(\bs{G_n}))_{j,k}|\leq t, \forall (j,k)\in\{1,...,M\}^2\}$. Under \ref{ass:RE}, on $\Omega_{G_n}$ , for all $J\subset\{1,...,M\}$ such that $|J|\leq s$ and for all $b\in \mathbb{R}^M\backslash\{0\}$ such that $||b_{J^c}||_1\leq a_0||b_J||_1$, we have
\begin{align*}
\bs{b^TG_nb}&=\bs{b}^T(\bs{G_n}-\mathbb{E}(\bs{G_n}))\bs{b}+\bs{b}^T\mathbb{E}(\bs{G_n})\bs{b}\\
&\geq \bs{b}^T(\bs{G_n}-\mathbb{E}(\bs{G_n}))\bs{b}+\bs{\kappa_0}^2||\bs{b}_J||_2^2.
\end{align*}
Since $\bs{b}^T(\bs{G_n}-\mathbb{E}(\bs{G_n}))\bs{b}=\sum_{j=1}^n\sum_{k=1}^n(\bs{G_n}-\mathbb{E}(\bs{G_n}))_{j,k}b_jb_k$, on $\Omega_{G_n}$, under \ref{ass:RE} we deduce that
$$\bs{b^TG_nb}\geq-\sum_{i,j}t|b_i||b_j|+\bs{\kappa_0}^2||\bs{b}_J||_2^2.$$
Since $||\bs{b}||_2\leq ||\bs{b}||_1\leq (1+a_0)||\bs{b}_J||_1\leq (1+a_0)\sqrt{s}||\bs{b}_J||_2$, we finally obtain
$$\bs{b^TG_nb}\geq (-t(1+a_0)^2s+\kappa_0^2)||\bs{b}_J||_2^2.$$
We choose $t=A_0\bs\kappa^2/(1+a_0)^2s$ with $\bs{\kappa}=\bs{\kappa_0}/\sqrt{2A_0}$ to get $\bs{b^TG_nb}\geq \bs{\kappa_0}^2||\bs{b}_J||_2^2$.

It remains to calculate $\mathbb{P}(\Omega_{G_n})$. The coefficient $(j,k)$ of the matrix $\bs{G_n}-\mathbb{E}(\bs{G_n})$ is given by
$$\dfrac{1}{n}\sum_{i=1}^n(\Lambda_i-\mathbb{E}(\Lambda_i))f_j(\bs{Z_i})f_k(\bs{Z_i}).$$
For sake of simplicity, we put $\zeta_i^{j,k}=\Lambda_if_j(\bs{Z_i})f_k(\bs{Z_i})$ for $i=1,...,n$ and $(j,k)\in\{1,...,M\}^2$ fixed.
To apply a standard Bernstein's inequality to the independent random variables $\zeta_1^{j,k},...,\zeta_n^{j,k}$, $(j,k)\in\{1,...,M\}^2$ we have to verify that 
$(1/n)\sum_{i=1}^{n}\mathbb{E}|\zeta_i^{j,k}|^m\leq m!vc^{m-2}$ for some positive constants $v$ and $c$ and for all integers $m\geq2$ (see Proposition 2.9 in Massart \cite{Massart}). Under Assumptions \ref{A_0} and \ref{ass:dico1} the variables $\zeta_i^{j,k}$ are bounded, $|\zeta_i^{j,k}|\leq A_0L^2$ for all $i=1,...,n$ and $(j,k)\in\{1,...,M\}^2$, so that the previous assumption is satisfied with $v=\sum_{i=1}^{n}\mathbb{E}[\zeta_i^{j,k}]\leq A_0^2L^4$ and $c=A_0L^2$ and the Bernstein's inequality applied to $(\zeta_i^{j,k})_{i=1,...,n}$ is
\begin{align}
\label{standard}
\mathbb{P}\Big(\sumin(\zeta_i^{j,k}-\mathbb{E}[\zeta_i^{j,k}])>x\Big)\leq \exp\Big(-\dfrac{x^2}{2(v+cx)}\Big).
\end{align}
From (\ref{standard}), we get
$$\mathbb{P}\Big(|(\bs{G_n}-\mathbb{E}(\bs{G_n}))_{i,j}|>\dfrac{A_0\bs{\kappa}^2}{(1+a_0)^2s}\Big)\leq 2\exp\Big(-\dfrac{n\bs{\kappa}^4}{2(1+a_0)^2sL^2(L^2(1+a_0)^2s+\bs{\kappa}^2/3)}\Big).$$
So the probability of $\Omega_{G_n}^c$ is given by
\begin{align*}
\mathbb{P}(\Omega_{G_n}^c)&=\mathbb{P}\Big(\exists (j,k)\in\{1,...,M\}^2 : |(\bs{G_n}-\mathbb{E}(\bs{G_n}))_{j,k}|>\dfrac{A_0\bs{\kappa}^2}{(1+a_0)^2s}\Big)\\
&\leq \displaystyle\sum_{j=1}^M\sum_{k=1}^M\mathbb{P}\Big(|(\bs{G_n}-\mathbb{E}(\bs{G_n}))_{j,k}|>\dfrac{A_0\bs{\kappa}^2}{(1+a_0)^2s}\Big)\\
&\leq 2M^2\exp\Big(-\dfrac{n\bs{\kappa}^4}{2(1+a_0)^2sL^2(L^2(1+a_0)^2s+\bs{\kappa}^2/3)}\Big),
\end{align*}
and by denoting 
$$\pi_n=2M^2\exp\Big(-\dfrac{n\bs{\kappa}^4}{2(1+a_0)^2sL^2(L^2(1+a_0)^2s+\bs{\kappa}^2/3)}\Big),$$
we finally get (\ref{eq:REn}) with probability larger than $1-\pi_n$. 
\qed

\subsection{Proof of Theorem \ref{thm:fast-oracle}}
Let introduce the event $\bs\Omega_{\bs{\RE_n}(s,a_0)}(\bs\kappa)=\Big\{0<\bs\kappa=\underset{|J|\leq s}{\underset{J\subset\{1,...,M\},}{\min}}\underset{||\bs{b}_{J^c}||_{1}\leq a_{0}||\bs{b}_{J}||_{1}}{\underset{\bs{b}\in\mathbb{R}^M\backslash\{0\},}{\min}}\dfrac{(\bs{b^TG_nb})^{1/2}}{||\bs{b}_{J}||_{2}}\Big\}$.
We start from Inequality (\ref{eq:B1}) and the fact that on $\mathcal{A}$, for $\bs\beta\in\Gamma(\mu)$,
$$|(\bs{\hat{\beta}_{L}^{\mu}}-\bs\beta)^T\bs{\eta_{n,\tau}}|\leq\sumip\dfrac{\omega_{j}}{2}|(\hat{\beta}_{L}^{\mu}-\beta)_{j}|.$$
It follows that
$$\widetilde{K}_{n}(\lambda_{0},\lambda_{\hat{\beta}_{L}^{\mu}})+\sumip\dfrac{\omega_{j}}{2}|(\hat{\beta}_{L}^{\mu}-\beta)_{j}|\leq\widetilde{K}_{n}(\lambda_{0},\lambda_{\beta})+\sumip\omega_{j}(|(\hat{\beta}_{L}^{\mu}-\beta)_{j}|+|\beta_{j}|-|(\hat{\beta}_{L}^{\mu})_{j}|).$$ 
On $J(\bs\beta)^c$, $|(\hat{\beta}_{L}^{\mu}-\beta)_{j}|+|\beta_{j}|-|(\hat{\beta}_{L}^{\mu})_{j}|=0$, so on $\mathcal{A}$ we obtain 
\begin{equation}
\label{eq:BRT1}
\widetilde{K}_{n}(\lambda_{0},\lambda_{\hat{\beta}_{L}^{\mu}})+\sumip\dfrac{\omega_{j}}{2}|(\hat{\beta}_{L}^{\mu}-\beta)_{j}|\leq\widetilde{K}_{n}(\lambda_{0},\lambda_{\beta})+2\sum_{j\in J(\beta)}\omega_{j}|(\hat{\beta}_{L}^{\mu}-\beta)_{j}|.
\end{equation}
We apply Cauchy-Schwarz Inequality to the second right hand side of (\ref{eq:BRT1}) to get 

\begin{equation}
\label{eq:BRT2}
\widetilde{K}_{n}(\lambda_{0},\lambda_{\hat{\beta}_{L}^{\mu}})+\sumip\dfrac{\omega_{j}}{2}|(\hat{\beta}_{L}^{\mu}-\beta)_{j}|\leq \widetilde{K}_{n}(\lambda_{0},\lambda_{\beta})+2\sqrt{|J(\bs\beta)|}\sqrt{\sum_{j\in J(\beta)}\omega_{j}^2|\hat\beta_{L}^{\mu}-\beta|_{j}^2}.
\end{equation}

With the notations $\bs\Delta=\bs D(\bs{\hat\beta_{L}^{\mu}}-\bs\beta)$ and $\bs D=(\diag(\omega_{j}))_{1\leq j\leq M}$ introduced in Subsection \ref{subsec:fast} , Inequalities (\ref{eq:BRT1}) and (\ref{eq:BRT2}) become 
\begin{equation}
\label{eq:BRT11}
\widetilde{K}_{n}(\lambda_{0},\lambda_{\hat{\beta}_{L}^{\mu}})+\dfrac{1}{2}||\bs\Delta||_{1}\leq\widetilde{K}_{n}(\lambda_{0},\lambda_{\beta})+2||\bs\Delta_{J(\beta)}||_{1},
\end{equation}
and 
\begin{equation}
\label{eq:BRT12}
\widetilde{K}_{n}(\lambda_{0},\lambda_{\hat{\beta}_{L}^{\mu}})\leq\widetilde{K}_{n}(\lambda_{0},\lambda_{\beta})+2\sqrt{|J(\beta)|}||\bs\Delta_{J(\beta)}||_{2}.
\end{equation}
Consider,
\begin{equation}
\label{eq:cas2}
\mathcal{A}_{1}=\{\zeta\widetilde{K}_{n}(\lambda_{0},\lambda_{\beta})\leq 2||\bs\Delta_{J(\beta)}||_{1}\}.
\end{equation}
On $\mathcal{A}\bigcap\mathcal{A}_{1}^c$, the result of the theorem follows immediately from (\ref{eq:BRT11}). 
As soon as, $||\bs\Delta_{J(\beta)^c}||_{1}\leq \left(3+4/\zeta\right)||\bs\Delta_{J(\beta)}||_{1},$
on $\bs\Omega_{\bs{\RE_n}(s,a_0)}(\bs\kappa)$, with $a_{0}=\left(3+4/\zeta\right)\text{ and } \bs{\kappa}=(1/\sqrt{2A_0})\bs\kappa_0(s,a_{0})$ we get 
$$\bs\kappa^2||\bs\Delta_{J(\bs\beta)}||_{2}^2\leq \bs\Delta^T\bs{G_{n}}\bs\Delta.$$
So, initially we will assume that $||\bs\Delta_{J(\beta)^c}||_{1}\leq \left(3+4/\zeta\right)||\bs\Delta_{J(\beta)}||_{1},$ and we will verify later that this inequality holds. Since,
\begin{align*}
\bs\Delta^T\bs{G_{n}}\bs\Delta&=\dfrac{1}{n}\sumin\Big(\displaystyle\sum_{j=1}^{M}\omega_j(\hat\beta_{L,j}^{\mu}-\beta_j)f_j(\bs{Z_{i}})\Big)^2\Lambda_{i}(\tau)\\
&\leq(\underset{1\leq j\leq M}{\max}\omega_j)^2\dfrac{1}{n} \sumin\inttau\Big(\log(\alpha_{0}(t)\e^{f_{\hat\beta_{L}^{\mu}}(\bs{Z_{i}})})-\log(\alpha_{0}(t)\e^{f_{\beta}(\bs{Z_{i}})})\Big)^2\diff \Lambda_{i}(t)\\
&\leq (\underset{1\leq j\leq M}{\max}\omega_j)^2||\log\lambda_{\hat\beta_{L}^{\mu}}-\log\lambda_{\beta}||_{n,\Lambda}^2,
\end{align*}
on $\bs\Omega_{\bs{\RE_n}(s,a_0)}(\bs\kappa)$, Inequality (\ref{eq:BRT12}) becomes on $\mathcal{A}\cap\bs\Omega_{\bs{\RE_n}(s,a_0)}(\bs\kappa)$
\begin{align*}
\widetilde{K}_{n}(\lambda_{0},\lambda_{\hat{\beta}_{L}^{\mu}})&\leq \widetilde{K}_{n}(\lambda_{0},\lambda_{\beta})+2\sqrt{|J(\bs\beta)|}(\underset{1\leq j\leq M}{\max}\omega_{j})\bs\kappa^{-1}||\log\lambda_{\hat\beta_{L}^{\mu}}-\log\lambda_{\beta}||_{n,\Lambda}\\
&\leq \widetilde{K}_{n}(\lambda_{0},\lambda_{\beta})+2\sqrt{|J(\bs\beta)|}(\underset{1\leq j\leq M}{\max}\omega_{j})\bs\kappa^{-1}(||\log\lambda_{\hat\beta_{L}^{\mu}}-\log\lambda_{0}||_{n,\Lambda}+||\log\lambda_{0}-\log\lambda_{\beta}||_{n,\Lambda}).
\end{align*}

The following proposition (proof in Annexe \ref{annexeA})  connects the weighted empirical norm and the empirical Kullback divergence. 
\begin{proposition}
\label{prop:comparaison}
Under Assumption \ref{ass:voisinage}, for all $\bs{\beta}\in\Gamma(\mu)$, 
\begin{equation*}
\mu'||\log{\lambda_{\beta}}-\log{\lambda_{0}}||^2_{n,\Lambda}\leq \widetilde{K}_{n}(\lambda_{0},\lambda_{\beta})\leq \mu''||\log{\lambda_{\beta}}-\log{\lambda_{0}}||^2_{n,\Lambda},
\end{equation*}
where $\mu'=\phi(\mu)/\mu^2$, $\mu''=\phi(-\mu)/\mu^2$ and $\phi(t)=\e^{-t}+t-1$.
\end{proposition}

Now, applying Proposition \ref{prop:comparaison}, it follows that
\begin{align*}
\widetilde{K}_{n}(\lambda_{0},\lambda_{\hat{\beta}_{L}^{\mu}})\leq\widetilde{K}_{n}(\lambda_{0},\lambda_{\beta})+2\sqrt{|J(\bs\beta)|}(\underset{1\leq j\leq M}{\max}\omega_{j})\dfrac{\bs\kappa^{-1}}{\sqrt{\mu'}}\left( \sqrt{\widetilde{K}_{n}(\lambda_{0},\lambda_{\hat\beta_{L}^{\mu}})}+\sqrt{\widetilde{K}_{n}(\lambda_{0},\lambda_{\beta})}\right).
\end{align*}
We now use the elementary inequality $2uv\leq bu^2+\dfrac{v^2}{b}$ with $b>1$, $u=\sqrt{|J(\bs\beta)|}(\underset{1\leq j\leq M}{\max}\omega_{j})\bs\kappa^{-1}$ and $v$ being either $\sqrt{\dfrac{1}{\mu'}\widetilde{K}_{n}(\lambda_{0},\lambda_{\hat\beta_{L}^{\mu}})}$ or $\sqrt{\dfrac{1}{\mu'}\widetilde{K}_{n}(\lambda_{0},\lambda_{\beta})}$. Consequently
$$\widetilde{K}_{n}(\lambda_{0},\lambda_{\hat{\beta}_{L}^{\mu}})\leq \widetilde{K}_{n}(\lambda_{0},\lambda_{\beta})+2b|J(\bs\beta)|(\underset{1\leq j\leq M}{\max}\omega_{j})^2\bs\kappa^{-2}+\dfrac{1}{b\mu'}\widetilde{K}_{n}(\lambda_{0},\lambda_{\hat\beta_{L}^{\mu}})+\dfrac{1}{b\mu'}\widetilde{K}_{n}(\lambda_{0},\lambda_{\beta}).$$
Hence,
$$\left(1-\dfrac{1}{\mu'b}\right)\widetilde{K}_{n}(\lambda_{0},\lambda_{\hat{\beta}_{L}^{\mu}})\leq \left(1+\dfrac{1}{b\mu'}\right)\widetilde{K}_{n}(\lambda_{0},\lambda_{\beta})+2b|J(\bs\beta)|(\underset{1\leq j\leq M}{\max}\omega_{j})^2\bs\kappa^{-2},$$
and 
$$\widetilde{K}_{n}(\lambda_{0},\lambda_{\hat{\beta}_{L}^{\mu}})\leq\dfrac{b\mu'+1}{b\mu'-1}\widetilde{K}_{n}(\lambda_{0},\lambda_{\beta})+2\dfrac{b^2\mu'}{b\mu'-1}|J(\beta)|(\underset{1\leq j\leq M}{\max}\omega_{j})^2\bs\kappa^{-2}.$$
We take $\dfrac{b\mu'+1}{b\mu'-1}=1+\zeta$ and $C(\zeta,\mu)=2\dfrac{b^2\mu'}{b\mu'+1}$ a constant depending on $\zeta$ and $\mu$. It follows that for any $\bs\beta\in\Gamma(\mu)$ : 
$$\widetilde{K}_{n}(\lambda_{0},\lambda_{\hat{\beta}_{L}^{\mu}})\leq (1+\zeta)\Big\{\widetilde{K}_{n}(\lambda_{0},\lambda_{\beta})+C(\zeta,\mu)|J(\bs\beta)|(\underset{1\leq j\leq M}{\max}\omega_{j})^2\bs\kappa^{-2}\Big\}.$$
Finally, taking the infimum over all $\bs\beta\in\Gamma(\mu)$ such that $|J(\bs\beta)|\leq s$, we obtain (\ref{eq:fast-oracle}).

We have now to verify that $||\bs\Delta_{J(\beta)^c}||_{1}\leq \left(3+4/\zeta\right)||\bs\Delta_{J(\beta)}||_{1}.$ On $\mathcal{A}\bigcap\mathcal{A}_{1}$, applying (\ref{eq:BRT11}) we get that
$$||\bs\Delta||_{1}\leq 4\left(1+\dfrac{1}{\zeta}\right)||\bs\Delta_{J(\beta)}||_{1},$$ so by splitting $\bs\Delta=\bs\Delta_{J(\bs\beta)}+\bs\Delta_{J(\bs\beta)^c}$, we finally obtain 
$$||\bs\Delta_{J(\beta)^c}||_{1}\leq \left(3+\dfrac{4}{\zeta}\right)||\bs\Delta_{J(\beta)}||_{1}.$$

Finally, Lemma \ref{REn} ensures that $\mathbb{P}(\mathcal{A}^c\cup\bs\Omega_{\bs{\RE_n}(s,a_0)}^c(\bs\kappa))\leq A_{\varepsilon,\nu}e^{-x} + \pi_n$, which achieves the proof of Theorem \ref{thm:fast-oracle}. 
\qed

\subsection{Proof of Corollary \ref{cor:corollary}}
Corollary \ref{cor:corollary} follows from Proposition \ref{prop:comparaison} and same arguments as in the proof of Theorem \ref{thm:fast-oracle} with
\[b=\dfrac{\mu'(1+\zeta)+\mu''}{\mu'(1+\zeta)-\mu''}.\]
\qed

\subsection{Proof of Theorem \ref{thm:selection}}
To prove Inequality (\ref{eq:prediction}) of Theorem \ref{thm:selection}, we start from (\ref{eq:BRT1}) with $\bs\beta=\bs\beta_{0}$ and $\bs{\hat\beta_L}$ defined by (\ref{eq:lassosimple}). Consequently $\widetilde{K}_{n}(\lambda_{0},\lambda_{\beta})=~0$.
Here we give the proposition that  gives the relation between the empirical Kullback divergence and the empirical norm, in the case of variable selection. 
\begin{proposition}
\label{relationsel}
Under Assumption (\ref{covborne}), there exist two positive numerical constants $\xi$ and $\xi'$ such that 
$$\xi ||\bs{(\hat\beta_L-\beta_0)^TX}||_{n,\Lambda}^2\leq \widetilde{K}_{n}(\lambda_0,\lambda_{\bs{\hat\beta_L}})\leq \xi'||\bs{(\hat\beta_L-\beta_0)^TX}||_{n,\Lambda}^2.$$
\end{proposition}
The proof of Proposition (\ref{relationsel}) is given in Annexe \ref{annexeB}. Applying Proposition \ref{relationsel} with $\lambda_{0}(t,\bs{Z_{i}})=\alpha_{0}(t)\e^{\bs{\beta_{0}^TZ_{i}}}$ and $\lambda_{\hat\beta_{L}}(t,\bs{Z_{i}})=\alpha_{0}(t)\e^{\bs{\hat\beta_{L}^TZ_{i}}}$, we obtain that, on $\mathcal{A}=\underset{j=1}{\overset{p}{\bigcap}}\left\{|\eta_{n,\tau}(f_{j})|\leq\Gamma_1\dfrac{\omega_{j}}{2}\right\}$
\begin{align}
\label{eq:selection1}
\xi||\bs{(\hat\beta_{L}-\beta_{0})^TX}||_{n,\Lambda}^2+\Gamma_1\sum_{j=1}^{p}\dfrac{\omega_{j}}{2}|\hat\beta_{L}-\beta_{0}|_{j}\leq 2\Gamma_1\sum_{j\in J_{0}}\omega_{j}|\hat\beta_{L}-\beta_{0}|_{j}.
\end{align}
From this inequality, we deduce 
\begin{align}
\label{eq:pred}
\xi||\bs {X(\hat\beta_L-\beta_0)}||_{n,\Lambda}^2\leq 2\Gamma_1\sum_{j\in J_{0}}\omega_{j}|\hat\beta_{L}-\beta_{0}|_{j}
                                                        \leq 2\sqrt{|J_{0}|}\Gamma_1||\bs{\Delta}_{0,J_{0}}||_{2}.
\end{align}
From (\ref{eq:selection1}), we also have
$$\displaystyle{\sum_{j=1}^{p}\omega_{j}|\hat\beta_{L}-\beta_{0}|_{j}}\leq4\sum_{j\in J_{0}}\omega_{j}|\hat\beta_{L}-\beta_{0}|_{j}$$
and we obtain
$||\bs{\Delta_{0}}||_{1}\leq 4||\bs{\Delta_{0}}_{J_{0}}||_{1}.$
We then split $||\bs{\Delta_{0}}||_{1}=||\bs{\Delta_{0}}_{J_{0}}||_{1}+||\bs{\Delta_{0}}_{J_{0}^c}||_{1}$ to get
\begin{align}
\label{eq:beforere}
||\bs{\Delta_{0}}_{J_{0}^c}||_{1}\leq3||\bs{\Delta_{0}}_{J_{0}}||_{1}.
\end{align}
On $\bs\Omega_{\bs{\RE_n}(s,a_0)}(\bs{\kappa'})$, with $a_0=3$ and $\bs{\kappa'}=(1/\sqrt{2A_0})\bs{\kappa_0}(s,3)$ we get 
\begin{align}
\label{eq:re}
||\bs {X\Delta_0}||_{n,\Lambda}^2\geq \bs{\kappa'}^2||\bs{\Delta}_{0,J_{0}}||_{2}^2.
\end{align}
According to (\ref{eq:pred}), we conclude that on $\mathcal{A}\cap \bs\Omega_{\bs{\RE_n}(s,a_0)}(\bs{\kappa'})$
$$\xi||\bs {X(\hat\beta_{L}-\beta_{0})}||_{n,\Lambda}^2\leq 2\sqrt{|J_{0}|}\Gamma_1\underset{1\leq j\leq p}{\max}\omega_{j}\dfrac{||\bs {X(\hat\beta_{L}-\beta_{0})}||_{n,\Lambda}}{\bs{\kappa'}},$$
which entails that
\begin{equation*}
||\bs X(\bs{\hat\beta_{L}}-\bs{\beta_{0}})||_{n,\Lambda}^2\leq \dfrac{4|J_{0}|}{\xi^2\bs{\kappa'}^2}\Gamma_1^2(\underset{1\leq j\leq p}{\max}\omega_{j})^2,
\end{equation*}
with $\mathbb{P}(\mathcal{A}\cap\bs\Omega_{\bs{\RE_n}(s,a_0)}(\bs\kappa'))\geq 1-A_{\varepsilon,\nu}e^{-\Gamma_1x}-\pi_n$.


Let us come to the proof of Inequality (\ref{eq:selection}) in Theorem \ref{thm:selection}. On $\mathcal{A}\cap\bs\Omega_{\bs{\RE_n}(s,a_0)}(\bs\kappa')$, with $a_0=3$, Inequality (\ref{eq:pred}) becomes 
$$\xi\dfrac{\bs{\kappa'}^2}{\underset{1\leq j\leq M}{\max}\omega_j^2}||\bs{\Delta}_{0,J_{0}}||_{2}^2\leq 2\sqrt{|J_{0}|}\Gamma_1||\bs{\Delta}_{0,J_{0}}||_{2},$$
and hence
\begin{align}
\label{eq:min}
||\bs{\Delta}_{0,J_{0}}||_{2}\leq \dfrac{2\sqrt{|J_{0}|}}{\xi\bs{\kappa'}^2}\Gamma_1\underset{1\leq j\leq p}{\max}\omega_{j}^2.
\end{align}
According to  (\ref{eq:beforere}) and thanks to Cauchy-Schwarz Inequality, we have
$$||\bs{\Delta_{0}}||_{1}=||\bs{\Delta_{0}}_{J_{0}}||_{1}+||\bs{\Delta_{0}}_{J_{0}^c}||_{1}\leq 4||\bs{\Delta_{0}}_{J_{0}}||_{1}\leq 4\sqrt{|J_{0}|}||\bs{\Delta_{0}}_{J_{0}}||_{2}.$$
From (\ref{eq:min}), we get 
$$\dfrac{||\bs{\Delta_{0}}||_{1}}{4\sqrt{|J_{0}|}}\leq \dfrac{2\sqrt{|J_{0}|}}{\xi\bs{\kappa'}^2}\Gamma_1\underset{1\leq j\leq p}{\max}\omega_{j}^2,$$
and finally
$$||\bs{\bs{\hat\beta_{L}-\beta_{0}}}||_{1}\leq 8\dfrac{|J_{0}|}{\xi\bs{\kappa'}^2}\Gamma_1\dfrac{\underset{1\leq j\leq p}{\max}\omega_j^2}{\underset{1\leq j\leq p}{\min}\omega_{j}},$$
with $\mathbb{P}(\mathcal{A}\cap\bs\Omega_{\bs{\RE_n}(s,a_0)}(\bs\kappa'))\geq 1-A_{\varepsilon,\nu}e^{-\Gamma_1x}-\pi_n$.
\qed

\subsection{Proof of Theorem \ref{thm:inconnu}}
The proof is very similar to the one of Theorem \ref{thm:slow-oracle}. We start from (\ref{eq:difference-risk}) and (\ref{eq:proctocontrol}), and write
\begin{align}
\nonumber
\widetilde{K}_{n}(\lambda_{0},\lambda_{\hat\beta_{L},\hat\gamma_{L}})\leq \widetilde{K}_{n}(\lambda_{0},\lambda_{\beta,\gamma}) &+(\bs{\hat\gamma_{L}}-\bs\gamma)^T\bs{\nu_{n,\tau}}+\pen(\bs\gamma)-\pen(\bs{\hat\gamma_{L}})\\
\label{eq:secondconnu}
&+(\bs{\hat\beta_{L}}-\bs\beta)^T\bs{\eta_{n,\tau}}+\pen(\bs\beta)-\pen(\bs{\hat\beta_{L}}).
\end{align}
Set $\mathcal{A}$ and $\mathcal{B}$ such that
\begin{equation}
\label{ensembles}
\mathcal{A}=\underset{j=1}{\overset{M}{\bigcap}}\Big\{|\eta_{n,\tau}(f_{j})|\leq\dfrac{\omega_{j}}{2}\Big\} \text{ and } \mathcal{B}=\underset{k=1}{\overset{N}{\bigcap}}\Big\{|\nu_{n,\tau}(\theta_{k})|\leq \dfrac{\delta_{k}}{2}\Big\}.
\end{equation}
We apply Theorem \ref{thm:empirical-bernstein} to bound up $\mathbb{P}(\mathcal{A}^c)$ and $\mathbb{P}(\mathcal{B}^c)$ and obtain that
$$\mathbb{P}(\mathcal{A}^c)\leq c_{3,\varepsilon,c_{\ell}}\e^{-x} \text{ and } \mathbb{P}(\mathcal{B}^c)\leq \tilde c_{3,\tilde\varepsilon,c_{\ell}'}\e^{-y}.$$
Hence for $A_{\varepsilon,\nu}=c_{3,\varepsilon,c_\ell} \text{ and } B_{\tilde\varepsilon,\tilde\nu}=\tilde c_{3,\tilde\varepsilon,c_\ell'}$, we have
\begin{equation}
\label{proba}
\mathbb{P}[(\mathcal{A}\cap\mathcal{B})^c]=\mathbb{P}(\mathcal{A}^c\cup\mathcal{B}^c)\leq \mathbb{P}(\mathcal{A}^c)+\mathbb{P}(\mathcal{B}^c)\leq  A_{\varepsilon,\nu}\e^{-x}+B_{\tilde\varepsilon,\tilde\nu}\e^{-y},
\end{equation}
with
$$A_{\varepsilon,\nu}=\dfrac{2}{\log(1+\varepsilon)}\log\Big(2+\dfrac{A_0(\nu/n+\Phi(\nu/n))}{x/n}\Big)+1  \hspace{0.1cm}\mbox{ and }\hspace{0.1cm} B_{\tilde\varepsilon,\tilde\nu}=\dfrac{2}{\log(1+\tilde\varepsilon)}\log\Big(2+\dfrac{A_0(\tilde\nu/n+\Phi(\tilde\nu/n))}{y/n}\Big)+1$$
On $\mathcal{A}\cap\mathcal{B}$ arguing as in the proof of Theorem \ref{thm:slow-oracle}, with probability larger than $1-A_{\varepsilon,\nu}\e^{-x}-B_{\tilde\varepsilon,\tilde\nu}\e^{-y}$, we finish the proof by writing (\ref{eq:slow-inconnu}).
\qed

\subsection{Proof of Theorem \ref{thm:fast-inconnu}}
Let introduce the event $\bs\Omega_{\bs{\widetilde\RE_n}(s,r_0)}(\bs{\tilde\kappa})=\Big\{0<\bs{\tilde\kappa}=\underset{|J|\leq s}{\underset{J\subset\{1,...,M\},}{\min}}\underset{||\bs{b}_{J^c}||_{1}\leq r_{0}||\bs{b}_{J}||_{1}}{\underset{\bs{b}\in\mathbb{R}^M\backslash\{0\},}{\min}}\dfrac{(\bs{b^T\tilde{G}_nb})^{1/2}}{||\bs{b}_{J}||_{2}}\Big\}$.
We start from Inequality (\ref{eq:secondconnu}). On $\mathcal{A}\cap\mathcal{B}$ defined in (\ref{ensembles}), for $(\bs\beta,\bs\gamma)\in\widetilde\Gamma(\rho)$,
$$|(\bs{\hat{\beta}_{L}}-\bs\beta)^T\bs{\eta_{n,\tau}}|\leq\sumip\dfrac{\omega_{j}}{2}|(\hat{\beta}_{L}-\beta)_{j}| \text{ and }|(\bs{\hat{\gamma}_{L}}-\bs\gamma)^T\bs{\nu_{n,\tau}}|\leq\sum_{k=1}^{N}\dfrac{\delta_{k}}{2}|(\hat{\gamma}_{L}-\gamma)_{k}|, $$
and therefore
\begin{align}
\nonumber
\widetilde{K}_{n}(\lambda_{0},\lambda_{\hat{\beta}_{L}^{\rho},\hat\gamma_{L}^{\rho}})&+\sumip\dfrac{\omega_{j}}{2}|(\hat{\beta}_{L}^{\rho}-\beta)_{j}|+\sum_{k=1}^{N}\dfrac{\delta_{k}}{2}|(\hat{\gamma}_{L}^{\rho}-\gamma)_{k}|\\
\label{eq:BRTinconnu1}
&\leq\widetilde{K}_{n}(\lambda_{0},\lambda_{\beta,\gamma})+2\sum_{j\in J(\bs\beta)}\omega_{j}|(\hat{\beta}_{L}^{\rho}-\beta)_{j}|+2\sum_{k\in J(\bs\gamma)}\delta_{k}|(\hat{\gamma}_{L}^{\rho}-\gamma)_{k}|.
\end{align}
We then apply Cauchy-Schwarz inequality to the second right-term of (\ref{eq:BRTinconnu1}) and obtain 
\begin{align}
\nonumber
\widetilde{K}_{n}(\lambda_{0},\lambda_{\hat{\beta}_{L}^{\rho},\hat\gamma_{L}^{\rho}})&+\sumip\dfrac{\omega_{j}}{2}|(\hat{\beta}_{L}^{\rho}-\beta)_{j}|+\sum_{k=1}^{N}\dfrac{\delta_{k}}{2}|(\hat{\gamma}_{L}^{\rho}-\gamma)_{k}|\\
\label{eq:BRTinconnu2}
&\leq \widetilde{K}_{n}(\lambda_{0},\lambda_{\beta,\gamma})+2\sqrt{|J(\bs\beta)|}\sqrt{\sum_{j\in J(\beta)}\omega_{j}^2|\hat\beta_{L}^{\rho}-\beta|_{j}^2}+2\sqrt{|J(\bs\gamma)|}\sqrt{\sum_{k\in J(\bs\gamma)}\delta_{k}^2|\hat\gamma_{L}^{\rho}-\gamma|_{k}^2}.
\end{align}
With the notation of Subsection \ref{subsec:fast-inconnu}, Inequality (\ref{eq:BRTinconnu1}) is rewritten as :
\begin{equation}
\label{eq:BRTinconnu11}
\widetilde{K}_{n}(\lambda_{0},\lambda_{\hat{\beta}_{L}^{\rho},\hat\gamma_{L}^{\rho}})+\dfrac{1}{2}||\bs{\tilde\Delta}||_{1}\leq\widetilde{K}_{n}(\lambda_{0},\lambda_{\beta,\gamma})+2||\bs{\tilde\Delta}_{J(\bs\beta),J(\bs\gamma)}||_{1},
\end{equation}
where $\bs{\tilde\Delta}_{J(\beta),J(\gamma)}=\bs {\tilde D}\begin{pmatrix}(\bs{\hat\beta_{L}^{\rho}}-\bs\beta)_{J(\bs\beta)}\\
(\bs{\hat\gamma_{L}^{\rho}}-\bs\gamma)_{J(\bs\gamma)}\end{pmatrix}$. In the same way, Inequality (\ref{eq:BRTinconnu2}) becomes :
\begin{equation}
\label{eq:BRTinconnu12}
\widetilde{K}_{n}(\lambda_{0},\lambda_{\hat{\beta}_{L}^{\rho},\hat\gamma_{L}^{\rho}})\leq\widetilde{K}_{n}(\lambda_{0},\lambda_{\beta,\gamma})+4\max\left(\sqrt{|J(\bs\beta)|},\sqrt{|J(\bs\gamma)|}\right)||\bs{\tilde\Delta}_{J(\bs\beta),J(\bs\gamma)}||_{2}.
\end{equation}
Consider
\begin{equation}
\label{eq:casinconnu2}
\mathcal{A}_{1}=\zeta\widetilde{K}_{n}(\lambda_{0},\lambda_{\beta,\gamma})\leq 2||\bs{\tilde\Delta}_{J(\bs\beta),J(\bs\gamma)}||_{1}.
\end{equation}
On $\mathcal{A}\cap\mathcal{B}\cap\mathcal{A}_{1}$, Inequality (\ref{eq:fast-inconnu}) in Theorem \ref{thm:fast-inconnu} follows immediately from (\ref{eq:BRTinconnu11}). 
As soon as, $||\bs{\tilde\Delta}_{J(\bs\beta)^c,J(\bs\gamma)^c}||_{1}\leq \left(3+8\max\left(\sqrt{|J(\bs\beta)|},\sqrt{|J(\bs\gamma)|}\right)/\zeta\right)||\bs{\tilde\Delta}_{J(\bs\beta),J(\bs\gamma)}||_{1},$ on $\bs\Omega_{\bs{\widetilde\RE_n}(s,r_0)}(\bs{\tilde\kappa})$, with 
$$\bs{\tilde\kappa}=(1/\sqrt{2})\bs{\tilde\kappa_0}(s,r_{0})\quad \text{and}\quad r_{0}=\left(3+8\max\left(\sqrt{|J(\bs\beta)|},\sqrt{|J(\bs\gamma)|}\right)/\zeta\right),$$
we get that
$$\bs{\tilde\kappa}^2||\bs{\tilde\Delta}_{J(\bs\beta),J(\bs\gamma)}||_{2}^2\leq \bs{\tilde\Delta^T\tilde G_{n}\tilde\Delta} \quad\text{with}\quad \bs{\tilde\Delta^T\tilde G_{n}\tilde\Delta}\leq \underset{1\leq k\leq N}{\underset{1\leq j\leq M}{\max}}\{\omega_{j},\delta_{k}\}||\log\lambda_{\hat\beta_{L}^{\rho},\hat\gamma_{L}^{\rho}}-\log\lambda_{\beta,\gamma}||_{n,\Lambda}^2.$$
On $\mathcal{A}\cap\mathcal{B}\cap\bs\Omega_{\bs{\widetilde\RE_n}(s,r_0)}(\bs{\tilde\kappa})$, Equation (\ref{eq:BRTinconnu12}) becomes
\begin{align*}
\widetilde{K}_{n}(\lambda_{0},\lambda_{\hat{\beta}_{L}^{\rho},\hat\gamma_{L}^{\rho}})&\leq\widetilde{K}_{n}(\lambda_{0},\lambda_{\beta,\gamma})+4\max\left(\sqrt{|J(\bs\beta)|},\sqrt{|J(\gamma)|}\right)||\bs{\tilde\Delta}_{J(\bs\beta),J(\bs\gamma)}||_{2}\\
&\leq \widetilde{K}_{n}(\lambda_{0},\lambda_{\beta})+4\max\left(\sqrt{|J(\bs\beta)|},\sqrt{|J(\bs\gamma)|}\right)\underset{1\leq k\leq N}{\underset{1\leq j\leq M}{\max}}\{\omega_{j},\delta_{k}\}\bs{\tilde\kappa}^{-1}||\log\lambda_{\hat\beta_{L}^{\rho},\hat\gamma_{L}^{\rho}}-\log\lambda_{\beta,\gamma}||_{n,\Lambda}.
\end{align*}
Using that $||\log\lambda_{\hat\beta_{L}^{\rho},\hat\gamma_{L}^{\rho}}-\log\lambda_{\beta,\gamma}||_{n,\Lambda}\leq||\log\lambda_{\hat\beta_{L}^{\rho},\hat\gamma_{L}^{\rho}}-\log\lambda_{0}||_{n,\Lambda}+||\log\lambda_{0}-\log\lambda_{\beta,\gamma}||_{n,\Lambda}$,
we obtain that $\widetilde{K}_{n}(\lambda_{0},\lambda_{\hat{\beta}_{L}^{\rho},\hat\gamma_{L}^{\rho}})$ is less than
 $$\widetilde{K}_{n}(\lambda_{0},\lambda_{\beta,\gamma})+4\max\left(\sqrt{|J(\bs\beta)|},\sqrt{|J(\bs\gamma)|}\right)\underset{1\leq k\leq N}{\underset{1\leq j\leq M}{\max}}\{\omega_{j},\delta_{k}\}\bs{\tilde\kappa}^{-1}(||\log\lambda_{\hat\beta_{L}^{\rho},\hat\gamma_{L}^{\rho}}-\log\lambda_{0}||_{n,\Lambda}+||\log\lambda_{0}-\log\lambda_{\beta,\gamma}||_{n,\Lambda}).$$
This inequality involves both oracle inequalities in empirical Kullback divergence and in weighted empirical norm. 
 
In the same way that Proposition \ref{prop:comparaison}, we obtain a Proposition that connect the empirical Kullback divergence and the weighted empirical norm.
\begin{proposition}
\label{prop:comparaison1}
Under Assumption \ref{ass:voisinage1}, for all $(\bs{\beta},\bs{\gamma})\in\widetilde\Gamma(\rho)$, 
\begin{equation*}
\rho'||\log{\lambda_{\beta,\gamma}}-\log{\lambda_{0}}||^2_{n,\Lambda}\leq \widetilde{K}_{n}(\lambda_{0},\lambda_{\beta,\gamma})\leq \rho''||\log{\lambda_{\beta,\gamma}}-\log{\lambda_{0}}||^2_{n,\Lambda},
\end{equation*}
where $\rho'=\phi(\rho)/\rho^2$, $\rho''=\phi(-\rho)/\rho^2$ and $\phi(t)=\e^{-t}+t-1$.
\end{proposition}

Applying Proposition \ref{prop:comparaison1}, we obtain that $\widetilde{K}_{n}(\lambda_{0},\lambda_{\hat{\beta}_{L}^{\rho},\hat\gamma_{L}^{\rho}})$ is less than
$$\widetilde{K}_{n}(\lambda_{0},\lambda_{\beta,\gamma})+4\max\left(\sqrt{|J(\bs\beta)|},\sqrt{|J(\bs\gamma)|}\right)\underset{1\leq k\leq N}{\underset{1\leq j\leq M}{\max}}\{\omega_{j},\delta_{k}\}\dfrac{\bs{\tilde\kappa}^{-1}}{\sqrt{\rho'}}\left( \sqrt{\widetilde{K}_{n}(\lambda_{0},\lambda_{\hat\beta_{L}^{\rho},\hat\gamma_{L}^{\rho}})}+\sqrt{\widetilde{K}_{n}(\lambda_{0},\lambda_{\beta,\gamma})}\right).$$
Using again $2uv\leq bu^2+\dfrac{v^2}{b}$ with $b>1$, $u=2\max\left(\sqrt{|J(\bs\beta)|},\sqrt{|J(\bs\gamma)|}\right)\underset{1\leq k\leq N}{\underset{1\leq j\leq M}{\max}}\{\omega_{j},\delta_{k}\}\bs{\tilde\kappa}^{-1}$ and $v$ being either $\sqrt{\dfrac{1}{\rho'}\widetilde{K}_{n}(\lambda_{0},\lambda_{\hat\beta_{L}^{\rho},\hat\gamma_{L}^{\rho}})}$ or $\sqrt{\dfrac{1}{\rho'}\widetilde{K}_{n}(\lambda_{0},\lambda_{\beta,\gamma})}$, we obtain 
\begin{align*}
\widetilde{K}_{n}(\lambda_{0},\lambda_{\hat{\beta}_{L}^{\rho},\hat\gamma_{L}^{\rho}})&\leq \widetilde{K}_{n}(\lambda_{0},\lambda_{\beta,\gamma})+8b\max(|J(\bs\beta)|,|J(\bs\gamma)|)\Big(\underset{1\leq k\leq N}{\underset{1\leq j\leq M}{\max}}\{\omega_{j},\delta_{k}\}\Big)^2\bs{\tilde\kappa}^{-2}\\
&+\dfrac{1}{b\rho'}\widetilde{K}_{n}(\lambda_{0},\lambda_{\hat\beta_{L}^{\rho},\hat\gamma_{L}^{\rho}})+\dfrac{1}{b\rho'}\widetilde{K}_{n}(\lambda_{0},\lambda_{\beta,\gamma}).
\end{align*}
Hence,
$$\left(1-\dfrac{1}{b\rho'}\right)\widetilde{K}_{n}(\lambda_{0},\lambda_{\hat{\beta}_{L}^{\rho},\hat\gamma_{L}^{\rho}})\leq \left(1+\dfrac{1}{b\rho'}\right)\widetilde{K}_{n}(\lambda_{0},\lambda_{\beta,\gamma})+8b\max(|J(\bs\beta)|,|J(\bs\gamma)|)\Big(\underset{1\leq k\leq N}{\underset{1\leq j\leq M}{\max}}\{\omega_{j},\delta_{k}\}\big)^2\bs{\tilde\kappa}^{-2},$$
and 
\begin{align}
\label{eq:coro}
\widetilde{K}_{n}(\lambda_{0},\lambda_{\hat{\beta}_{L}^{\rho},\hat\gamma_{L}^{\rho}})\leq\dfrac{b\rho'+1}{b\rho'-1}\widetilde{K}_{n}(\lambda_{0},\lambda_{\beta,\gamma})+8\dfrac{b^2\rho'}{b\rho'-1}\max(|J(\bs\beta)|,|J(\bs\gamma)|)\Big(\underset{1\leq k\leq N}{\underset{1\leq j\leq M}{\max}}\{\omega_{j},\delta_{k}\}\Big)^2\dfrac{\bs{\tilde\kappa}^{-2}}{\rho'}.
\end{align}
We take $\dfrac{b\rho'+1}{b\rho'-1}=1+\zeta$ and we introduce $\widetilde{C}(\zeta,\rho)=8\dfrac{b^2\rho'}{b\rho'+1}$ a constant depending on $\zeta$ and $\rho$. For all $(\bs\beta,\bs\gamma)$ in $\widetilde\Gamma(\rho)$, we obtain
$$\widetilde{K}_{n}(\lambda_{0},\lambda_{\hat{\beta}_{L}^{\rho},\hat\gamma_{L}^{\rho}})\leq (1+\zeta)\Big\{\widetilde{K}_{n}(\lambda_{0},\lambda_{\beta,\gamma})+\widetilde{C}(\zeta,\rho)\max(|J(\bs\beta)|,|J(\gamma)|)\Big(\underset{1\leq k\leq N}{\underset{1\leq j\leq M}{\max}}\{\omega_{j},\delta_{k}\}\Big)^2\bs{\tilde\kappa}^{-2}\Big\}.$$
Finally, taking the infimum over all $(\bs\beta,\bs\gamma)\in\widetilde\Gamma(\rho)$ such that $\max(|J(\bs\beta)|,|J(\bs\gamma)|)\leq s$, we obtain Inequality (\ref{eq:fast-inconnu}).
Inequality (\ref{eq:fast-inconnu-emp}) follows by applying Proposition \ref{prop:comparaison} with $b=\dfrac{(1+\zeta)\rho'+\rho''}{(1+\zeta)\rho'-\rho''}$ in (\ref{eq:coro}).

We have now to verify that $||\bs{\tilde\Delta}_{J(\bs\beta)^c,J(\bs\gamma)^c}||_{1}\leq \left(3+8\max\left(\sqrt{|J(\bs\beta)|},\sqrt{|J(\bs\gamma)|}\right)/\zeta\right)||\bs{\tilde\Delta}_{J(\bs\beta),J(\bs\gamma)}||_{1}.$ We deduce from (\ref{eq:BRTinconnu11}) that, on $\mathcal{A}\cap\mathcal{B}\cap\mathcal{A}_{1}$,
$$||\bs{\tilde\Delta}||_{1}\leq 4\left(1+\dfrac{2}{\zeta}\max\left(\sqrt{|J(\bs\beta)|},\sqrt{|J(\bs\gamma)|}\right)\right)||\bs{\tilde\Delta}_{J(\bs\beta),J(\bs\gamma)}||_{1}.$$
By splitting $\bs{\tilde\Delta}=\bs{\tilde\Delta}_{J(\bs\beta),J(\bs\gamma)}+\bs{\tilde\Delta}_{J(\bs\beta)^c,J(\bs\gamma)^c}$, we infer that 
$$||\bs{\tilde\Delta}_{J(\bs\beta)^c,J(\bs\gamma)^c}||_{1}\leq \left(3+\dfrac{8}{\zeta}\max\left(\sqrt{|J(\bs\beta)|},\sqrt{|J(\bs\gamma)|}\right)\right)||\bs{\tilde\Delta}_{J(\bs\beta),J(\bs\gamma)}||_{1}.$$

To achieve the proof of Theorem \ref{thm:fast-inconnu}, we combine Equation (\ref{proba}) with Lemma \ref{tildeREn} to conclude
$$\mathbb{P}\Big[\Big(\mathcal{A}\cap\mathcal{B}\cap\bs\Omega_{\bs{\widetilde\RE_n}(s,r_0)}(\bs{\tilde\kappa})\Big)^c\Big]\leq A_{\varepsilon,\nu}\e^{-x}+B_{\tilde\varepsilon,\tilde\nu}\e^{-y}+\tilde\pi_n.$$
\qed

\subsection{Proof of Theorem \ref{thm:empirical-bernstein}}

The proofs of (\ref{eq:empirical-bernstein}) and (\ref{eq:empirical-bernstein1}) are quite similar, so we only present the one of (\ref{eq:empirical-bernstein}). To prove (\ref{eq:empirical-bernstein1}), it suffices to replace $\eta_{n,t}(f_{j})$ by the process $\nu_{n,t}(\theta_{k})$ throughout the following.
Denote by $U_{n,t}$ and $H_{i}(f_{j})$ the quantities
$$U_{n,t}(f_{j})=\dfrac{1}{n}\sumin\intt H_{i}(f_{j})\diff M_{i}(s) \text{ and } H_{i}(f_{j}):=\dfrac{f_{j}(\bs{Z_{i}})}{\underset{1\leq i\leq n}{\max}|f_{j}(\bs{Z_{i}})|}.$$ 
Since $H_{i}(f_{j})$ is a bounded predictable process with respect to $\mathcal{F}_{t}$, $U_{n,t}(f_{j})$ is a square integrable martingale. Its predictable variation is given by
$$\vartheta_{n,t}(f_{j})=n<U_{n}(f_{j})>_{t}=\dfrac{1}{n}\sumin\intt(H_i(f_j))^2\diff \Lambda_{i}(s)$$
and the optional variation of $U_{n,t}(f_{j})$ is
$$\hat\vartheta_{n,t}(f_{j})=n[U_{n}(f_{j})]_{t}=\dfrac{1}{n}\sumin\intt(H_i(f_j))^2\diff N_{i}(s).$$
We also define 
\begin{align}
\label{decalage}
\mathcal{\hat W}_n^{\nu}(f_j)=\dfrac{\nu/n}{\nu/n-\Phi(\nu/n)}\hat\vartheta_{n,t}(f_{j})+\dfrac{x/n}{\nu/n-\Phi(\nu/n)},
\end{align}
for $\nu\in(0,3)$ such that $\nu>\Phi(\nu)$ with $\Phi(u)=\e^{u}-u-1$.

From Inequality (7.12) in Hansen et al. \cite{Hansen}, for any $0<v<\omega<+\infty$, we have
\begin{align}
\label{7.12}
\mathbb{P}\Big(U_{n,t}(f_j)\geq \sqrt{\dfrac{2(1+\varepsilon)\mathcal{\hat W}_n^{\nu}(f_j)x}{n}}+\dfrac{x}{3n}, v\leq \mathcal{\hat W}_n^{\nu}(f_j)\leq \omega\Big)\leq 2\Big(\dfrac{\log(\omega/v)}{\log(1+\varepsilon)}+1\Big)\e^{-x}.
\end{align}
We focus now on removing the event $\{v\leq \mathcal{\hat W}_n^{\nu}(f_j)\leq \omega\}$ in (\ref{7.12}).
Let us consider the martingale given $\mathcal{F}_t$
\begin{align*}
\hat\vartheta_{n,t}(f_{j})-\vartheta_{n,t}(f_{j})&=\dfrac{1}{n}\sumin\intt(H_{i}(f_{j}))^2\Big(\diff N_{i}(s)-\diff\Lambda_{i}(s)\Big)=\dfrac{1}{n}\sumin\intt(H_{i}(f_{j}))^2\diff M_{i}(s),
\end{align*}
and let 
$$S_{\nu,t}(f_j)=\sumin\intt\Phi\big(\dfrac{\nu}{n}H_i^2(f_j)\Big)\diff \Lambda_i(s).$$
From van de Geer \cite{SVG}, we know that
$$\exp(\nu(\hat\vartheta_{n,t}(f_{j})-\vartheta_{n,t}(f_{j}))-S_{\nu,t}(f_j))$$
is a supermartingale. Now from Markov Inequality, for any $\nu,x>0$, we obtain that
\begin{align}
\label{S_nu}
\mathbb{P}\Big[|\hat\vartheta_{n,t}(f_{j})-\vartheta_{n,t}(f_{j})|\geq \dfrac{S_{\nu,t}(f_j)}{\nu}+\dfrac{x}{n}\Big]\leq 2\e^{-x}.
\end{align}
For any $0<h<1$ and $x>0$, $\Phi(xh)\leq h^2\Phi(x)$. This combined with the fact that $0<H^2_i(f_j)<1$, we get
\begin{align}
\label{Phi}
S_{\nu,t}(f_j)\leq \Phi(\nu/n)\sumin\intt H_i^4(f_j)\diff M_i(s)\leq \Phi(\nu/n)n\vartheta_{n,t}(f_j).
\end{align}
Combining (\ref{S_nu}) and (\ref{Phi}), we deduce that
\begin{align}
\label{DiffMart}
\mathbb{P}\Big[|\hat\vartheta_{n,t}(f_{j})-\vartheta_{n,t}(f_{j})|\geq \dfrac{\Phi(\nu/n)}{\nu/n}\vartheta_{n,t}(f_{j})+\dfrac{x}{\nu}\Big]\leq 2\e^{-x}.
\end{align}

Now, under Assumption \ref{A_0}, we have $\vartheta_{n,t}(f_{j})\leq A_0$, so the events 
$$\Omega_n^{\nu}=\Big\{\dfrac{x/n}{\nu/n-\Phi(\nu/n)}\leq \mathcal{\hat W}_n^{\nu}(f_j)\Big\}\cap\{\vartheta_{n,t}(f_{j})\leq A_0\}$$
is of probability one and thus
\begin{align}
\mathbb{P}\Big(U_{n,t}(f_j)\geq \sqrt{\dfrac{2(1+\varepsilon)\mathcal{\hat W}_n^{\nu}(f_j)x}{n}}+\dfrac{x}{3n}\Big)\leq\mathbb{P}\Big(\Big\{U_{n,t}(f_j)\geq \sqrt{\dfrac{2(1+\varepsilon)\mathcal{\hat W}_n^{\nu}(f_j)x}{n}}+\dfrac{x}{3n}\Big\} \cap \Omega_n^{\nu} \Big).
\end{align}
From (\ref{DiffMart}), we have 
$$\mathbb{P}\Big[\hat\vartheta_{n,t}(f_{j})\geq \vartheta_{n,t}(f_{j})\Big(1+\dfrac{\Phi(\nu/n)}{\nu/n}\Big)+\dfrac{x}{\nu}\Big]\leq \e^{-x},$$
and if we denote $E_n^{\nu}$ the event
$$E_n^{\nu}=\Big\{\hat\vartheta_{n,t}(f_{j})\leq \vartheta_{n,t}(f_{j})\Big(1+\dfrac{\Phi(\nu/n)}{\nu/n}\Big)+\dfrac{x}{\nu}\Big\},$$
we get 
\begin{align*}
\mathbb{P}\Big[U_{n,t}(f_j)\geq \sqrt{\dfrac{2(1+\varepsilon)\mathcal{\hat W}_n^{\nu}(f_j)x}{n}}+\dfrac{x}{3n}\Big]\leq \e^{-x}+\mathbb{P}\Big[\Big\{U_{n,t}(f_j)\geq \sqrt{\dfrac{2(1+\varepsilon)\mathcal{\hat W}_n^{\nu}(f_j)x}{n}}+\dfrac{x}{3n}\Big\} \cap \Omega_n^{\nu} \cap E_n^{\nu}\Big].
\end{align*}
On the event $E_n^{\nu}\cap\Omega_n^{\nu}$, from the definition of $\mathcal{\hat W}_n^{\nu}(f_j)$ given by (\ref{decalage}), we have
\begin{align}
\nonumber
\mathcal{\hat W}_n^{\nu}(f_j)&\leq \dfrac{\nu/n}{\nu/n-\Phi(\nu/n)}\Big(\vartheta_{n,t}(f_{j})\Big(1+\dfrac{\Phi(\nu/n)}{\nu/n}\Big)+\dfrac{x}{\nu}\Big)+\dfrac{x/n}{\nu/n-\Phi(\nu/n)}\\
\label{EetO}
&\leq A_0\dfrac{\nu/n+\Phi(\nu/n)}{\nu/n-\Phi(\nu/n)}+2\dfrac{x/n}{\nu/n-\Phi(\nu/n)}. 
\end{align}
From (\ref{EetO}), we obtain
\begin{align*}
&\mathbb{P}\Big[\Big\{U_{n,t}(f_j)\geq \sqrt{\dfrac{2(1+\varepsilon)\mathcal{\hat W}_n^{\nu}(f_j)x}{n}}+\dfrac{x}{3n}\Big\} \cap\Omega_n^{\nu} \cap E_n^{\nu}\Big]\\
\leq&\mathbb{P}\Big[U_{n,t}(f_j)\geq \sqrt{\dfrac{2(1+\varepsilon)\mathcal{\hat W}_n^{\nu}(f_j)x}{n}}+\dfrac{x}{3n}, \hspace{0.2cm}\dfrac{x/n}{\nu/n-\Phi(\nu/n)}\leq \mathcal{\hat W}_n^{\nu}(f_j)\leq A_0\dfrac{\nu/n+\Phi(\nu/n)}{\nu/n-\Phi(\nu/n)}+2\dfrac{x/n}{\nu/n-\Phi(\nu/n)}\Big].
\end{align*}
We now apply Inequality (\ref{7.12}) with $v=\dfrac{x/n}{\nu/n-\Phi(\nu/n)}$ and $\omega= A_0\dfrac{\nu/n+\Phi(\nu/n)}{\nu/n-\Phi(\nu/n)}+2\dfrac{x/n}{\nu/n-\Phi(\nu/n)}$,
\begin{align} 
\mathbb{P}\Big[U_{n,t}(f_j)\geq \sqrt{\dfrac{2(1+\varepsilon)\mathcal{\hat W}_n^{\nu}(f_j)x}{n}}+\dfrac{x}{3n}\Big]&\leq \e^{-x}+2\Big(\log\Big(\dfrac{A_0\dfrac{\nu/n+\Phi(\nu/n)}{\nu/n-\Phi(\nu/n)}+2\dfrac{x/n}{\nu/n-\Phi(\nu/n)}}{\dfrac{x/n}{\nu/n-\Phi(\nu/n)}}\Big)+1\Big)\e^{-x},\\
&\leq \Big(\dfrac{2}{\log(1+\varepsilon)}\log\Big(2+\dfrac{A_0(\nu/n+\Phi(\nu/n))}{x/n}\Big)+1\Big)\e^{-x}.
\end{align}
Now it suffices to multiply both sides of the inequality inside the probability by $||f_j||_{n,\infty}=\underset{1\leq i\leq n}{\max}|f_j(\bs{Z_i})|$ to end up the proof of Theorem \ref{thm:empirical-bernstein}.
\qed

\appendix

\section{Proof of Proposition \ref{prop:comparaison1}}
\label{annexeA}

The proof of Proposition \ref{prop:comparaison} and Proposition \ref{prop:comparaison1} are similar. So we only dprove Proposition \ref{prop:comparaison1} which corresponds to the general case.
To compare the empirical Kullback divergence (\ref{eq:Kullback}) and the weighted empirical norm (\ref{eq:norm}), we use Lemma 1 in Bach \cite{Bach}, that we recall here :

\begin{lemma}
\label{lem:taylor}
Let $g$ be a convex three times differentiable function $g:\mathbb{R}\rightarrow\mathbb{R}$ such that for all $t\in\mathbb{R}$, $|g'''(t)|\leq Sg''(t)$, for some $S\geq 0$. Then, for all $t\geq 0$ :
\begin{equation*}
\dfrac{g''(0)}{S^2}\phi(St)\leq g(t)-g(0)-g'(0)t\leq \dfrac{g''(0)}{S^2}\phi(-St) \text{ with } \phi(u)=\e^{-u}+u-1
\end{equation*}
\end{lemma}

This Lemma gives upper and lower Taylor expansions for some convex and three times differentiable function. It has been introduced to extend tools from self-concordant functions (i.e. which verify $|g'''(t)|\leq 2g''(t)^{3/2}$) and provide simple extensions of theoretical results for the square loss for logistic regression.

Let $h$ be a function on $[0,\tau]\times\mathbb{R}^p$ and define
$$G(h)=-\dfrac{1}{n}\sumin\int_{0}^{\tau}h(s,\bs{Z_{i}})\diff \Lambda_{i}(s)+\dfrac{1}{n}\sumin \inttau \e^{h(s,\bs{Z_{i}})}Y_{i}(s)\diff s.$$
Consider the function $g:\mathbb{R}\rightarrow\mathbb{R}$ defined by $g(t)=G(h+tk)$, where $h$ and $k$ are two functions defined on $\mathbb{R}^p$. By differentiating $G$ with respect to $t$ we get :
\begin{align*}
g'(t)&=-\dfrac{1}{n}\sumin\int_{0}^{\tau}k(s,\bs{Z_{i}})\diff \Lambda_{i}(s)+\dfrac{1}{n}\sumin \inttau k(s,\bs{Z_{i}})\e^{h(s,\bs{Z_{i}})+tk(s,\bs{Z_{i}})}Y_{i}(s)\diff s,\\
g''(t)&=\dfrac{1}{n}\sumin \inttau(k(s,\bs{Z_{i}}))^2\e^{h(s,\bs{Z_{i}})+tk(s,\bs{Z_{i}})}Y_{i}(s)\diff s,\\
g'''(t)&=\dfrac{1}{n}\sumin \inttau(k(s,\bs{Z_{i}}))^3\e^{h(s,\bs{Z_{i}})+tk(s,\bs{Z_{i}})}Y_{i}(s)\diff s.
\end{align*}
It follows that
$$|g'''(t)|\leq||k||_{n,\infty} g''(t).$$
Applying Lemma \ref{lem:taylor} with $S=||k||_{n,\infty}$, we obtain for all $t\geq 0$,
$$\dfrac{g''(0)}{||k||_{n,\infty}^2}\phi(t||k||_{n,\infty})\leq g(t)-g(0)-g'(0)t\leq \dfrac{g''(0)}{||k||_{n,\infty}^2}\phi(-t||k||_{n,\infty}).$$
Take $t=1$, $h(s,\bs{Z_{i}})=\log\lambda_{0}(s,\bs{Z_{i}}) \text{ and  for } (\bs\beta,\bs\gamma)\in\widetilde\Gamma(\rho), k(s,\bs{Z_{i}})=\log\lambda_{\beta,\gamma}(s,\bs{Z_{i}})-\log\lambda_{0}(s,\bs{Z_{i}}).$
We obtain 
\begin{equation}
\label{eq:bach1}
g''(0)\dfrac{\phi(||\log\lambda_{\beta,\gamma}-\log\lambda_{0}||_{n,\infty})}{||\log\lambda_{\beta,\gamma}-\log\lambda_{0}||_{n,\infty}^2}\leq G(\log\lambda_{\beta,\gamma})- G(\log\lambda_{0})-g'(0)\leq g''(0)\dfrac{\phi(-||\log\lambda_{\beta,\gamma}-\log\lambda_{0}||_{n,\infty})}{||\log\lambda_{\beta,\gamma}-\log\lambda_{0}||_{n,\infty}^2}.
\end{equation}
Now straightforward calculations show that g'(0)=0 and
\begin{align*}
g''(0)&=\dfrac{1}{n}\sumin \inttau((\log\lambda_{\beta,\gamma}-\log\lambda_{0})(s,\bs{Z_{i}}))^2\diff \Lambda_{i}(s)\\
&=||\log\lambda_{\beta,\gamma}-\log\lambda_{0}||_{n,\Lambda}^2.
\end{align*}
Replacing $g'(0)$ and $g''(0)$ by their expressions in (\ref{eq:bach1}) and noting that  
$$G(\log\lambda_{\beta,\gamma})-G(\log\lambda_{0})=\widetilde{K}_{n}(\lambda_{0},\lambda_{\beta,\gamma}),$$ 
we get

$$\dfrac{\phi(||\log\lambda_{\beta,\gamma}-\log\lambda_{0}||_{n,\infty})}{||\log\lambda_{\beta,\gamma}-\log\lambda_{0}||_{n,\infty}^2}||\log\lambda_{\beta,\gamma}-\log\lambda_{0}||_{n,\Lambda}^2\leq \widetilde{K}_{n}(\lambda_{0},\lambda_{\beta,\gamma})\leq \dfrac{\phi(-||\log\lambda_{\beta,\gamma}-\log\lambda_{0}||_{n,\infty})}{||\log\lambda_{\beta,\gamma}-\log\lambda_{0}||_{n,\infty}^2}||\log\lambda_{\beta,\gamma}-\log\lambda_{0}||_{n,\Lambda}^2.$$
According to Assumption \ref{ass:voisinage} for $(\bs\beta,\bs\gamma)\in\widetilde\Gamma(\rho)$, 
$$||\log\lambda_{\beta,\gamma}-\log\lambda_{0}||_{n,\infty}\leq \rho.$$
Since $\phi(t)/t^2$ is decreasing and bounded below by 0, we can deduce that
$$\dfrac{\phi(||\log\lambda_{\beta,\gamma}-\log\lambda_{0}||_{n,\infty})}{||\log\lambda_{\beta,\gamma}-\log\lambda_{0}||_{n,\infty}^2}\geq \dfrac{\phi(\rho)}{\rho^2}$$
and
$$\dfrac{\phi(-||\log\lambda_{\beta,\gamma}-\log\lambda_{0}||_{n,\infty})}{||\log\lambda_{\beta,\gamma}-\log\lambda_{0}||_{n,\infty}^2}\leq \dfrac{\phi(-\rho)}{\rho^2}.$$
Take $\rho':=\phi(\rho)/\rho^2>0$ and $\rho'':=\phi(-\rho)/\rho^2>0$ to finish the proof.
\qed

\section{Proof of Proposition \ref{relationsel}}
\label{annexeB}
The beginning of this proof is similar to the proof of Proposition \ref{prop:comparaison1}. 

$\bullet$ For $\bs\beta$ and $\bs\eta$ in $\mathbb{R}^M$, let $G: \mathbb{R}^M\rightarrow \mathbb{R}$ and $g: \mathbb{R}\rightarrow\mathbb{R}$ define by
$$G(\bs{\beta})=-\dfrac{1}{n}\displaystyle\sum_{i=1}^{n}\int_0^{\tau}\log(\alpha_0(s)\e^{\bs{\beta^TZ_i}})\diff \Lambda_i(s)+\dfrac{1}{n}\displaystyle\sum_{i=1}^{n}\int_{0}^{\tau}\alpha_0(s)\e^{\bs{\beta^TZ_i}}Y_i(s)\diff s \text{ and }g(t)=G(\bs\beta+t\bs\eta).$$
 By differentiating $G$ with respect to $t$, we get 
\begin{align*}
g'(t)&=-\dfrac{1}{n}\displaystyle\sum_{i=1}^{n}\int_{0}^{\tau}\bs{\eta^T\bs{Z_i}}\diff\Lambda_i(s)+\dfrac{1}{n}\displaystyle\sum_{i=1}^{n}\int_{0}^{\tau}\alpha_{0}(s)\bs{\eta^TZ_i}\e^{(\bs\beta+t\bs\eta)^T\bs{Z_i}}Y_i(s)\diff s\\
g''(t)&=\dfrac{1}{n}\displaystyle\sum_{i=1}^{n}\int_{0}^{\tau}\alpha_0(s)(\bs{\eta^TZ_i})^2\e^{(\bs\beta+t\bs\eta)^T\bs{Z_i}}Y_i(s)\diff s\\
g'''(t)&=\dfrac{1}{n}\displaystyle\sum_{i=1}^{n}\int_{0}^{\tau}\alpha_0(s)(\bs{\eta^TZ_i})^3\e^{(\bs\beta+t\bs\eta)^T\bs{Z_i}}Y_i(s)\diff s
\end{align*}
It follows that 
$$|g'''(t)|\leq \dfrac{1}{n}\displaystyle\sum_{i=1}^{n}\int_{0}^{\tau}\alpha_0(s)||\bs\eta||_2||\bs{Z_i}||_2(\bs{\eta^TZ_i})^2\e^{(\bs\beta+t\bs\eta)^T\bs{Z_i}}Y_i(s)\diff s,$$

Under Assumption \ref{covborne}, we can deduce that $|g'''(t)|\leq R||\bs\eta||_2g''(t).$ 
Now applying Lemma \ref{lem:taylor} with $S=R||\bs\eta||_2$, we obtain for all $t\geq 0$,
$$\dfrac{g''(0)}{R^2||\bs\eta||_2^2}\phi(R||\bs\eta||_2t)\leq g(t)-g(0)-g'(0)t\leq \dfrac{g''(0)}{R^2||\bs\eta||_2^2}\phi(-R||\bs\eta||_2t)$$
Take $t=1$, $\bs\beta=\bs{\beta_0}$ and $\bs\eta=\bs{\hat\beta_L-\beta_0}$,
to write
\begin{align}
\label{Bach1}
g''(0)\dfrac{\phi(R||\bs{\hat\beta_L-\beta_0}||_2)}{R^2||\bs{\hat\beta_L-\beta_0}||_2^2}\leq G(\bs{\hat\beta_L})-G(\bs{\beta_0})-g'(0)\leq g''(0)\dfrac{\phi(-R||\bs{\hat\beta_L-\beta_0}||_2)}{R^2||\bs{\hat\beta_L-\beta_0}||_2^2}
\end{align}
Now straightforward calculations show that $g'(0)=0$ and
\begin{align*}
g''(0)=\dfrac{1}{n}\displaystyle\sum_{i=1}^{n}\int_0^{\tau}((\bs{\hat\beta_L}-\bs{\beta_0})^T\bs{Z_i})^2\alpha_0(s)\e^{\bs{\beta_0^TZ_i}}Y_i(t)\diff t=\dfrac{1}{n}\displaystyle\sum_{i=1}^{n}\int_0^{\tau}((\bs{\hat\beta_L}-\bs{\beta_0})^T\bs{Z_i})^2\diff \Lambda_i(s)=||(\bs{\hat\beta_L}-\bs{\beta_0})^T\bs{X}||_{n,\Lambda}^2,
\end{align*}
Replacing $g'(0)$ and $g''(0)$ by their expressions in (\ref{Bach1}) and noting that 
$$G(\bs{\hat\beta_L})-G(\bs{\beta_0})=\widetilde{K}_n(\lambda_0,\lambda_{\bs{\hat\beta_L}}),$$
we get 
\begin{equation}
\label{INKNE}
||(\bs{\hat\beta_L}-\bs{\beta_0})^T\bs{X}||_{n,\Lambda}^2\dfrac{\phi(R||\bs{\hat\beta_L-\beta_0}||_2)}{R^2||\bs{\hat\beta_L-\beta_0}||_2^2}\leq \widetilde{K}_n(\lambda_0,\lambda_{\bs{\hat\beta_L}})\leq ||(\bs{\hat\beta_L}-\bs{\beta_0})^T\bs{X}||_{n,\Lambda}^2\dfrac{\phi(-R||\bs{\hat\beta_L-\beta_0}||_2)}{R^2||\bs{\hat\beta_L-\beta_0}||_2^2}.
\end{equation}

$\bullet$ Now, we will show that $R||\bs{\hat\beta_L-\beta_0}||_2$ is bounded.
From Equation (\ref{eq:BRT11}) with $\bs{\hat\beta_L^{\mu}}=\bs{\hat\beta_L}$ and $\bs\beta=\bs{\beta_0}$, we can deduce that 
$$\widetilde{K}_{n}(\lambda_0,\lambda_{\bs{\hat\beta}_L})\leq \dfrac{3}{2}\Gamma_1||\bs\Delta_0||_1,$$
where $\bs{\Delta_0}=\bs{D(\hat\beta_L-\beta_0)}$ and $\bs D=(\diag(\omega_j))_{1\leq j\leq M}$.
From (\ref{INKNE}), we have
\begin{align*}
\widetilde{K}_{n}(\lambda_0,\lambda_{\bs{\hat\beta}_L})&\geq \dfrac{||\bs{(\hat\beta_L-\beta_0)^TX}||_{n,\Lambda}^2}{R^2||\bs{(\hat\beta_L-\beta_0)}||_2^2}\phi(R||\bs{(\hat\beta_L-\beta_0)}||_2)
\end{align*}
We apply Assumption \ref{ass:RE} with $a_0=3$ and $\bs{\kappa'}=\bs\kappa'(s,3)$
and we infer that 
$$\bs\kappa'^2||\bs{\Delta}_{0,J_0}||_{2}^2\leq ||\bs{\Delta_0^TX}||_{n,\Lambda}^2.$$
So we have,
$$\dfrac{\bs\kappa'^2||\bs{\Delta}_{0,J_0}||_2^2}{\underset{1\leq j\leq M}{\max}\omega_j^2}\dfrac{\phi(R||\bs{\hat\beta_L-\beta_0}||_2)}{R^2||\bs{\hat\beta_L-\beta_0}||_2^2}\leq  \dfrac{3}{2}\Gamma_1||\bs{\Delta_0}||_1.$$
We can now use, with $s=|J_0|$, $||\bs{\Delta_0}||_2\leq ||\bs{\Delta_0}||_1\leq 4||\bs{\Delta}_{0,J_0}||_1\leq 4\sqrt{s}||\bs{\Delta}_{0,J_0}||_2$ to get
\begin{align*}
\bs\kappa'^2\phi(R||\bs{\hat\beta_L-\beta_0}||_2)&\leq \dfrac{3}{2}\Gamma_1\dfrac{\underset{1\leq j\leq M}{\max}\omega_j^2}{\underset{1\leq j\leq M}{\min}\omega_j^2}\underset{1\leq j\leq M}{\max}\omega_j\dfrac{(4\sqrt{s}||(\bs{\hat\beta_L-\beta_0})_{J_0}||_2)^2R^2||\bs{\hat\beta_L-\beta_0}||_2}{||(\bs{\hat\beta_L-\beta_0})_{J_0}||_2^2}\\
&\leq 24\Gamma_1\dfrac{\underset{1\leq j\leq M}{\max}\omega_j^2}{\underset{1\leq j\leq M}{\min}\omega_j^2}\underset{1\leq j\leq M}{\max} \omega_j sR^2||\bs{\Delta_0}||_2.
\end{align*}
A short calculation shows that for all $k\in(0,1]$ : 
$$\e^{-2k(1-k)^{-1}}+(1-k)2k(1-k)^{-1}-1\geq 0.$$
(see Bach \cite{Bach} for more details)
So by taking $2k(1-k)^{-1}=R||\bs{\hat\beta_L-\beta_0}||_2$, we have
$$\e^{-R||\bs{\hat\beta_L-\beta_0}||_2}+R||\bs{\hat\beta_L-\beta_0}||_2-1\geq \dfrac{R^2||\bs{\hat\beta_L-\beta_0}||_2^2}{2+R||\bs{\hat\beta_L-\beta_0}||_2}$$
and we deduce that
$$\dfrac{\bs\kappa'^2R^2||\bs{\hat\beta_L-\beta_0}||_2^2}{2+R||\bs{\hat\beta_L-\beta_0}||_2}\leq 24\Gamma_1\dfrac{\underset{1\leq j\leq M}{\max}\omega_j^2}{\underset{1\leq j\leq M}{\min}\omega_j^2}\underset{1\leq j\leq M}{\max} \omega_jsR^2||\bs{\hat\beta_L-\beta_0}||_2.$$
This implies that $R||\bs{\hat\beta_L-\beta_0}||_2\leq \dfrac{\dfrac{48\Gamma_1Rs}{\bs\kappa'^2}\dfrac{\underset{1\leq j\leq M}{\max}\omega_j^2}{\underset{1\leq j\leq M}{\min}\omega_j^2}\underset{1\leq j\leq M}{\max}\omega_j^2}{1-\dfrac{24\Gamma_1Rs}{\bs\kappa'^2}\dfrac{\underset{1\leq j\leq M}{\max}\omega_j^2}{\underset{1\leq j\leq M}{\min}\omega_j^2}\underset{1\leq j\leq M}{\max}\omega_j^2}\leq 2$ as soon as $\Gamma_1\leq \dfrac{1}{48Rs}\dfrac{\underset{1\leq j\leq M}{\min}\omega_j^2}{\underset{1\leq j\leq M}{\max}\omega_j^2}\dfrac{\bs\kappa'^{2}}{\underset{1\leq j\leq M}{\max}\omega_j}$. 

$\bullet$ Since $\phi(t)/t^2$ is decreasing and bounded below by $0$, we can deduce that 
$$\dfrac{\phi(R||\bs{\hat\beta_L-\beta_0}||_2)}{R^2||\bs{\hat\beta_L-\beta_0}||_2^2}\geq\dfrac{\phi(2)}{4}$$
and
$$\dfrac{\phi(-R||\bs{\hat\beta_L-\beta_0}||_2)}{R^2||\bs{\hat\beta_L-\beta_0}||_2^2}\leq\dfrac{\phi(-2)}{4}$$
Take $\xi:=\phi(2)/4>0$ and $\xi':=\phi(-2)/4>0$ and conclude that
$$\xi||\bs{(\hat\beta_L-\beta_0)^TX}||_{n,\Lambda}^2\leq \widetilde{K}_n(\lambda_0,\lambda_{\bs{\hat\beta_L}})\leq\xi'||\bs{(\hat\beta_L-\beta_0)^TX}||_{n,\Lambda}^2.$$
\qed
%
%
%
\newpage
%
%
%
%
%
\bibliography{biblio2}
\bibliographystyle{plain}
\end{document}